\documentclass{tran-l}
\usepackage{latexsym,color}
\newtheorem{theorem}{Theorem}
\newtheorem{lemma}[theorem]{Lemma}
\newtheorem{corollary}[theorem]{Corollary}
\newcommand{\bWn}{{\bf W}^M}
\newcommand{\bWinf}{{\bf W}^\infty}
\newcommand{\bVn}{{\bf V}^M}
\newcommand{\bVinf}{{\bf V}^\infty}

\newcommand{\eqdef}{\, =\kern -13pt\raise 6pt\hbox{{\tiny\textrm{def}}}\,\,}
\begin{document}

\title{Closed form summation of $C$-finite sequences}
\author{ Curtis Greene}
\address{Haverford College, Haverford, PA 19041-1392}
\email{cgreene@haverford.edu}
\author{Herbert S. Wilf}
\address{University of Pennsylvania, Philadelphia, PA 19104-6395}
\email{wilf@math.upenn.edu}
\date{May 15, 2004, and, in revised form, December 9, 2004.}
\dedicatory{To David Robbins}

\begin{abstract}
 We consider sums of
the form
\[\sum_{j=0}^{n-1}F_1(a_1n+b_1j+c_1)F_2(a_2n+b_2j+c_2)\dots
F_k(a_kn+b_kj+c_k),\] in which each $\{F_i(n)\}$ is a sequence that satisfies a linear
recurrence of degree $D(i)<\infty$, with constant coefficients. We assume further that the
$a_i$'s and the $a_i+b_i$'s are
all nonnegative integers. We prove that such a sum always has a closed
form, in the sense that it evaluates to a linear combination of a finite set of
monomials in the values of the sequences $\{F_i(n)\}$ with coefficients
that are polynomials in $n$. We
explicitly describe two different sets of monomials that will form
such a linear combination, and give an algorithm for finding these closed
forms, thereby completely automating the solution of this class of
summation problems. We exhibit tools for determining when these
explicit evaluations are unique of their type, and prove that in a
number of interesting cases they are indeed unique. We also discuss
some special features of the case of ``indefinite summation," in which $a_1=a_2=\cdots = a_k = 0$.
\end{abstract}

\maketitle

\tableofcontents

\section{Introduction}
In section 1.6 of \cite{pwz} the following assertion is made:
\begin{quotation}
{\small All Fibonacci number identities such as Cassini's
$F_{n+1}F_{n-1}-F_n^2=(-1)^n$
(and much more complicated ones), are routinely provable
using Binet's formula\index{Binet's formula}:
\[
F_n:={{1} \over {\sqrt{5}}} \left(
\left({{1+ \sqrt{5}} \over {2}}\right)^n-
\left({{1-\sqrt{5}} \over {2}}\right)^n \right).
\]}
\end{quotation}
This is followed by a brief Maple program that proves Cassini's
identity by
substituting Binet's formula on the left side and showing that it then
reduces to $(-1)^n$. Another method of proving these identities is
given in \cite{z1}, in which it is observed that one can find the recurrence
relations that are satisfied by each of the two sides of the identity
in question, show that they are the same and that the initial values
agree, and the identity will then be proved.

The purpose of this note is to elaborate on these ideas by showing
how to \textit{derive}, instead of only to verify, summation
identities for a certain class of sequence sums, and to show that
this class of sums always has closed form in a certain sense, and
that these closed forms can be found entirely algorithmically.
Indeed, a Mathematica program that carries out the procedures that
we develop in this paper can be downloaded from the web sites of
the authors \cite{gwmma}.

We deal with the class of $C$-finite sequences (see \cite{z1}). These
are the sequences $\{F(n)\}_{n\ge 0}$ that satisfy linear recurrences of
fixed span with constant coefficients. The Fibonacci numbers,
e.g., will do nicely for a prototype sequence of this kind.
The kind of sum that we will consider first will be of the form
(\ref{eq:sumdef}) below. We will say that such a sum has an {\it
$F$-closed form} if there is a linear combination of a fixed
(i.e., independent of $n$)  number  of monomials in values of the
$F$'s such that for all $n$ the sum $f(n)$ is equal to that linear
combination.

For example, consider the sum
$$
f(n) = \sum_{j=0}^{n-1} F(j)^2 F(2n-j)
$$
where the $F$'s are Fibonacci numbers.
In Section 3.1 we will see how to use our method to show that $f(n)$
can be expressed in the form (\ref{eq:coll2}), which
is a linear combination of five monomials in the $F$'s. Hence the sum $f(n)$
has an $F$-closed form.

More generally, we will consider functions $F(n)$ satisfying a
recurrence of minimal order $D$ with constant coefficients, whose
associated polynomial has roots $r_1,r_2,\dots,r_d$, of
multiplicities $e_1,e_2,\dots,e_d$, where $\sum_i e_i=D$. Such a
function may be expressed in the form
\begin{equation}\label{eq:fform}
F(n) = \sum_{m=1}^{d} \sum_{h=0}^{e_{m}-1}\,\lambda_{m,h}\,n^h \,(r_m)^n,
\end{equation}
where the $r_i$ are distinct and nonzero, and
$\lambda_{m,e_m-1}\not=0$ for all $m$.

We will begin by evaluating sums of the form
\begin{equation}
f(n)=\sum_{j=0}^{n-1}F(a_1n+b_1j+c_1)\cdots
F(a_kn+b_kj+c_k)\label{eq:sumdef}
\end{equation}
in which the $a$'s, $b$'s, and $c$'s are given integers. We assume
further that, for all $i$, $a_i\ge 0$ and $a_i+b_i\ge 0$
and at least one of these is positive. Later we
will generalize this result to allow the $F$'s in the summand to
be different $C$-finite functions.
The principal result of this paper is perhaps Theorem \ref{th:mpl} below, which
proves  in full generality, i.e., with arbitrary root multiplicities and
with the $F$'s in the summand being all different, the existence of closed
forms, and exhibits an explicit finite basis for the solution space.

It is elementary and well known
that $f(n)$ is $C$-finite, and one can readily obtain explicit
expressions for $f(n)$ in terms of the roots $r_m$. Our first
results show how to obtain formul\ae\ for sums $f(n)$ of the form
(\ref{eq:sumdef}) as a polynomial
in the $F$'s, based on two different explicit sets of ``target''
monomials in the $F$'s. Using the first target set, we obtain the
following result.

\begin{theorem}\label{thm1}
The sum $f(n)$ in (\ref{eq:sumdef}) has an $F$-closed form. It is
equal to a linear combination of
monomials in the $F$'s, of the form
\begin{eqnarray}
&F((a_1+b_1)n+i_1)\dots F((a_k+b_k)n+i_k),&\qquad 0\le i_{\nu}\le
d-1\nonumber\\
&\psi_{i_1,\dots,i_k}(n)F(a_1n+i_1)F(a_2n+i_2)\dots F(a_kn+i_k),&\qquad 0\le
i_{\nu}\le d-1,\label{eq:ans1}
\end{eqnarray}
where $\psi_{i_1,\dots,i_k}(n)$ denotes a polynomial of degree at
most $\beta = 1 +  \Delta \mu$, where $\Delta = \max_m\{e_m-1\}$
and $\mu =|\{i: a_i=0\}|$. If $F$ is rational-valued then there
are solutions in which all coefficients to be determined are rational.
\end{theorem}
We note that, if the roots $r_i$ are distinct, then $\beta = 1$ and
the polynomials $\psi_{i_1,\dots,i_k}(n)$ are linear.
The next theorem gives a closed form in terms of an
alternate target set of monomials.
\begin{theorem}\label{thm2}
The sum $f(n)$ in (\ref{eq:sumdef}) can be expressed in $F$-closed form as a
linear combination
of monomials of the form
\begin{eqnarray}
&F(n+i_1)F(n+i_2)\dots F(n+i_P),&\qquad 0\le
i_{\nu}\le d-1,\nonumber\\
&\psi_{i_1,\dots,i_k}(n)F(n+i_1)F(n+i_2)\dots F(n+i_Q),&\qquad 0\le
i_{\nu}\le d-1,\label{eq:ans3}
\end{eqnarray}
where
$P=(a_1+b_1)+(a_2+b_2)+\cdots+(a_k+b_k)$, $Q=a_1+a_2+\cdots+ a_k$, and
$\psi_{i_1,\dots,i_k}(n)$ denotes a polynomial of degree
$\gamma = \max\{0,1+ \Delta (k-\sum\{a_i | a_i>0\})\}$.
If $F$ is rational-valued then there are solutions in which all coefficients to be determined are rational.
\end{theorem}
For example, when $F(n)$ is the $n$th Fibonacci number,
Theorem 2 states that any sum of the form (\ref{eq:sumdef}) can
be expressed as a linear combination of monomials in $F(n)$ and
$F(n+1)$, with rational linear polynomial coefficients, where
those monomials have at most two distinct degrees. Again, we note that, if the roots $r_i$ are distinct, then
$\gamma=1$ and the polynomials $\psi_{i_1,\dots,i_k}(n)$ are
linear.

The natural domain for these questions is the vector space
$\bVinf$ of complex-valued functions on $\{0,1,2,\dots \}$.
However, to obtain our expansions it is only necessary to work in
the vector space $\bVn$ of functions on $\{1,\dots,M\}$, where $M$
is the number of unknown coefficients to be determined.
More
precisely, we define $M$ to be equal to the number of
``algebraically distinct'' monomials of the form $n^h \prod_\nu
F((a_\nu+b_\nu)n + i_\nu)$ or $n^h \prod_\nu F(a_\nu n + i_\nu)$
generated by (\ref{eq:ans1}) or (\ref{eq:ans3}). Here we consider
two monomials to be equivalent if they differ by a rearrangement
of factors, or by a constant factor arising from cases where
$a_\nu + b_\nu=0$ or $a_\nu=0$. (For examples, see Section 3.)
Then we have the following result.

\begin{theorem}\label{lem3}
Let $\bWn\subseteq \bVn$ be the vector space of complex-valued functions on
$\{1,2,\dots,M \}$ spanned by the monomials in (\ref{eq:ans1}),
where $M$ is the number of algebraically distinct monomials generated by
(\ref{eq:ans1}), as defined
in the previous paragraph. Let $\bWinf\subseteq \bVinf$ be the vector space of functions
on $\{0,1,\dots\}$ spanned by the same monomials. If two linear
combinations of monomials of type (\ref{eq:ans1}) agree in $\bWn$,
then they agree in $\bWinf$. A similar statement holds for
monomials of type (\ref{eq:ans3}).
\end{theorem}
As a consequence, we can obtain expressions of type (\ref{eq:ans1}) or (\ref{eq:ans3}) by equating
$M$ values of $f(n)$ to the values of the assumed linear combinations, and solving for the
coefficients. We note that $M\le (\beta+2)d^k$ in case (\ref{eq:ans1}) and $M\le d^P + (\gamma+1)d^Q$
in case (\ref{eq:ans3}).

In general, $F$-closed expressions are not unique. For example, we
may add terms of the form $\Psi(F) (F(n+2)- F(n+1) - F(n))$, where
$\Psi(F)$ is any polynomial in the $F(a n+i)$, to an expression
involving Fibonacci numbers and get another valid $F$-closed form.
However, the formats described by (\ref{eq:ans1}) and
(\ref{eq:ans3}) are highly restrictive, and the
resulting expressions can be shown to be unique in a surprising
number of cases. We will return to the question of uniqueness and,
more generally, to the problem of computing $\dim(\bWinf)$, in
Sections 4-6.

The structure of this paper is as follows.
In Sections 2-6 we consider
the case where the roots $r_i$ all have multiplicity one. In this
case, both the statements and proofs of our results are
considerably simpler,
and will serve as models for the more general case to
be presented later. Section 2 gives proofs of
Theorems \ref{thm1}, \ref{thm2} and \ref{lem3} in this special case.
Section 3 illustrates these results with several examples. Sections 4-6
consider the issue of uniqueness and dimension, again for the distinct root case.
Section 7 drops the assumption of distinct roots, and gives a proof of
Theorems \ref{thm1} and \ref{thm2} in a more general form
(Theorem \ref{th:mpl}) where the
factors in the
summand of (\ref{eq:sumdef}) may involve different $F$'s. Section 8 considers
some issues that arise when $a_1=a_2=\cdots=a_k=0$, i.e., when the problem of
computing $f(n)$ is an indefinite summation problem. Section 9 contains
more examples.

\section{Proofs in the case of distinct roots}
In this section, we will assume that $F(n)$ satisfies a recurrence of minimal order
$d$, with
distinct roots, and hence can be expressed in the form
\begin{equation}\label{eq:fform2}
F(n) = \sum_{m=1}^d \lambda_m r_m^n
\end{equation}
with the $r_m$ distinct and the $\lambda_m$ nonzero.
Expanding the right side of (\ref{eq:sumdef}) above and using
(\ref{eq:fform2}), we find that
\[f(n)=\sum_{j=0}^{n-1}\prod_{\ell=1}^k\left\{\sum_{m=1}^d
\lambda_mr_m^{a_{\ell}n+b_{\ell}j+c_{\ell}}  \right\}.\]
A typical term in the expansion of the product will look like
\begin{equation}\label{eq:typi}
Kr_{m_1}^{a_1n+b_1j+c_1}r_{m_2}^{a_2n+b_2j+c_2}\dots
r_{m_k}^{a_kn+b_kj+c_k},
\end{equation}
in which $K$ is a constant, i.e., is independent of $n$ and $j$, which
may
be different at different places in the exposition below.
Since we are about to sum the above over $j=0\dots n-1$, put
\[\Theta=r_{m_1}^{b_1}r_{m_2}^{b_2}\dots r_{m_k}^{b_k},\]
because this is the quantity that is raised to the $j$th power in the
expression (\ref{eq:typi}). Now there are two cases, namely $\Theta=1$
and
$\Theta\neq 1$.

Suppose $\Theta=1$. Then the sum of our typical term (\ref{eq:typi})
over
$j=0\dots n-1$ is
\begin{equation}\label{eq:theta1}Kn\left(r_{m_1}^{a_1}r_{m_2}^{a_2}\dots
r_{m_k}^{a_k}\right)^n.\end{equation}
On the other hand, if $\Theta\neq 1$ then the sum of our typical term
(\ref{eq:typi}) over $j=0\dots n-1$ is
\begin{equation}\label{eq:thetaneq1}
K\left\{\left(r_{m_1}^{a_1+b_1}\dots
r_{m_k}^{a_k+b_k}\right)^n-\left(r_{m_1}^{a_1}\dots
r_{m_k}^{a_k}\right)^n\right\}.
\end{equation}

The next task will be to express these results in terms of various
members
of the sequence $\{F(n)\}$ instead of in terms of various powers of the
$r_i$'s. To  do that we write out (\ref{eq:fform2}) for $d$ consecutive
values of $n$, getting
\begin{eqnarray*}
F(n+i)&=&\sum_{m=1}^d\lambda_mr_m^{n+i}\qquad (i=0,1,\dots, d-1)\\
&=&\sum_{m=1}^d(\lambda_mr_m^i)r_m^n\qquad(i=0,1,\dots, d-1).
\end{eqnarray*}
We regard these as $d$ simultaneous linear equations in the unknowns
$\{r_1^n,\dots,r_d^n\}$, with a coefficient matrix that is a
nonsingular
diagonal matrix times a Vandermonde based on distinct points, and is
therefore nonsingular. Hence for each $m=1,\dots, d$, $r_m^n$
is a
linear combination of $F(n),F(n+1),\dots,F(n+d-1)$, with coefficients
that
are independent of $n$. Thus in eqs. (\ref{eq:theta1}),
(\ref{eq:thetaneq1}) we can replace each $r_{m_i}^{na_i}$ by a linear
combination of $F(a_in),F(a_in+1),\dots,F(a_in+d-1)$ and we can replace
each
$r_{m_i}^{n(a_i+b_i)}$ by a linear combination of
$F((a_i+b_i)n),F((a_i+b_i)n+1),\dots,F((a_i+b_i)n+d-1)$.

After making these replacements, we see that the two possible
expressions
(\ref{eq:theta1}), (\ref{eq:thetaneq1}) contribute monomials that
are all of the form (\ref{eq:ans1}), with the polynomials
$\psi_{i_1,\dots,i_k}$ all linear.
This establishes the existence of expansions in monomials of type (\ref{eq:ans1}),
as claimed in Theorem 1.

To prove the corresponding claim made in Theorem \ref{thm2},
it suffices to observe that, in the above argument, we
could
have written $r_{m_i}^{na_i}= (r_{m_i}^n)^{a_i}$ and replaced it by a
homogeneous
polynomial of degree $a_i$ in $F(n),F(n+1),\dots,F(n+d-1)$.  Similar
reasoning applies to
$r_{m_i}^{n(a_i+b_i)}$. Thus all of the resulting monomials are of type
(\ref{eq:ans3}).

We continue now with the proof of Theorem \ref{lem3}. The arguments are identical
for cases (\ref{eq:ans1}) and (\ref{eq:ans3}), so we will consider
only case (\ref{eq:ans1}).
We have observed that for each $i=0,\dots,d-1$, $F(n + i)$
is a linear combination of $r_1^n,\dots,r_d^n$, and
conversely. Hence, in both $\bVn$ and $\bVinf$, the linear span of the set
$$F(a_1n+i_1)F(a_2n+i_2)\dots F(a_kn+i_k), \qquad 0\le
i_{\nu}\le d-1$$
is equal to the linear span of the set
$\{\theta_1^n,\theta_2^n,\dots,\theta_{d^k}^n\}$, where
the $\theta_j$ range over all monomials of the form
$$r^{a_1}_{m_1}r^{a_2}_{m_2}\cdots r^{a_k}_{m_k}.$$
Arguing similarly for the other cases,
we see that the linear span of all monomials of type (\ref{eq:ans1}) is
equal to the linear span of the set of $3d^k$ functions
\begin{eqnarray}
&\theta_1^n,\theta_2^n,\dots,\theta_{d^k}^n&\nonumber\\
&\psi_1^n,\psi_2^n,\dots,\psi_{d^k}^n&\label{eq:wpq1}\\
&n\psi_1^n,n\psi_2^n,\dots,n\psi_{d^k}^n,& \nonumber
\end{eqnarray}
where the $\theta_i$ are as defined above and the $\psi_j$
range over all monomials of the form
$$r^{a_1+b_1}_{m_1}r^{a_2+b_2}_{m_2}\cdots r^{a_k+b_k}_{m_k}.$$
We claim
that the number of {\it distinct} functions appearing in (\ref{eq:wpq1}) is
less than or equal to $M$. Indeed, it is straightforward to check that the map
$F(\theta n + i)\mapsto (r_{i+1}^\theta)^n$ extends to a well-defined, surjective map from the set of
equivalence classes of
monomials of type (\ref{eq:ans1}) to the set of functions appearing in (\ref{eq:wpq1}).

Now suppose that $\Phi(n)$ and $\Psi(n)$ are linear combinations of monomials of type
(\ref{eq:ans1}), with
$\Phi(n)=\Psi(n)$ for $n=1,2,\dots,M$. We know that
$\Phi(n)$ and $\Psi(n)$ can both be expressed in the form
$$
\sum_{i} c_i\, \theta_i^n + \sum_{j} d_j \,\psi_j^n + \sum_{k} e_k\, n \psi_k^n
$$
for some constants $c_i, d_j, e_k$, where the sum is over distinct elements of
(\ref{eq:wpq1}) and hence there are at most $M$ terms in the sum. It follows from standard
results in the theory of difference equations (e.g., see \cite{jordan}, Chapter 11) that
$\Phi(n)$ and $\Psi(n)$ satisfy the same linear recurrence of order at most $M$ with
constant coefficients, e.g., the recurrence with characteristic polynomial
$\prod_i (t-\theta_i) \prod_j (t-\psi_j)^2$. Hence the values
of $\Phi(n)$ and $\Psi(n)$ are
completely determined by their values for $n=1,2,\dots,M$, and since
they agree for these values, they must agree for all $n$. This completes
the proof of Theorem \ref{lem3}. $\Box$

The proof also shows that the C-finite degree of $f(n)$ is bounded
by $3d^k$. Sharper bounds appear in Corollaries \ref{gendegbound}
and \ref{gendegboundx} below.

\section{Examples}
\subsection{A Fibonacci sum}
This work was started when a colleague asked about the sum
\begin{equation}\label{eq:coll}
f(n)=\sum_{j=0}^{n-1}F(j)^2F(2n-j),\end{equation}
in which the $F$'s are the Fibonacci numbers. If we refer to the
general
form (\ref{eq:sumdef}) of the question we see that in this case
\[k=3;\,d=2;\,(a_1,b_1,c_1)=(a_2,b_2,c_2)=(0,1,0);\,(a_3,b_3,c_3)=(2,-1,0).\]
If we now refer to the general form
(\ref{eq:ans1}) of the
answer we see that the sum $f(n)$ is a linear combination of monomials
\[nF(2n),\,F(2n),\,nF(2n+1),\,F(2n+1),\,F(n)^3,\,F(n)^2F(n+1),\,F(n)F(n+1)^2,\,F(n+1)^3.\]
Hence we assume a linear combination of these monomials and equate its
values to those of $f(n)$ for $n=0,1,\dots,7$ to determine the
constants of
the linear combination. The result is that
\begin{equation}\label{eq:coll2}
f(n)=\frac12\left(F(2n)+F(n)^2F(n+1)-F(n)F(n+1)^2+F(n+1)^3-F(2n+1)\right).
\end{equation}
This formula is expressed in terms of
monomials of type (\ref{eq:ans1}). Using monomials of type (\ref{eq:ans3}),
we obtain the alternate expression
\begin{eqnarray}
f(n) &=& \frac12 \left(2 F(n)F(n+1) - 2 F(n)^2 - F(n+1)^2 + F(n)^2 F(n+1)
\right. \label{eq:coll2a}\\
&&\qquad\qquad\qquad\qquad\qquad\qquad\qquad \left.- F(n)F(n+1)^2+
F(n+1)^3 \right)\nonumber
\end{eqnarray}
In Section 5 we will show that both of these expression are unique, i.e.,
(\ref{eq:coll2})
is the unique $F$-closed formula for $f(n)$ of type (\ref{eq:ans1}) and
(\ref{eq:coll2a}) is the unique $F$-closed
formula of type (\ref{eq:ans3}).

\subsection{An example involving subword avoidance}
Given an alphabet of $A\ge 2$ letters,  let $W$ be some fixed word of
three
letters such that no proper suffix of $W$ is also a proper prefix of
$W$.
For example, $W=aab$ will do nicely. Let $G(n)$ be the number of
$n$-letter
words over $A$ that do not contain $W$ as a subword. It is well known,
and
obvious, that
\begin{equation}\label{eq:subwrecur}
G(n)=AG(n-1)-G(n-3),
\end{equation}
with $G(0)=1, G(1)=A, G(2) = A^2$,
so this is a $C$-finite sequence. It is easy to check
that the roots of its associated polynomial equation are distinct
for all $A\geq 2$.
Suppose we want to evaluate the sum $g(n)=\sum_{j=0}^{n-1}G(j)^2.$ Then
\[k=2;\,d=3;\,(a_1,b_1,c_1)=(a_2,b_2,c_2)=(0,1,0).\]
Using either Theorem 1 or Theorem 2, we see that $g(n)$ is a linear
combination of the monomials
\[1,\,n,\,G(n)^2,\,G(n)G(n+1),\,G(n)G(n+2),\,G(n+1)G(n+2),\,G(n+1)^2,\,G(n+2)^2.\]
As before, we assume a linear combination of these monomials with
constants
to be determined, and we equate the result to
computed values of $g(n)$, for
$n=0,1,\dots,7$,
to solve for the constants. The end result is that

\begin{eqnarray}
g(n)&=&\frac{1}{A(A-2)}\biggl(1-(A-1)^2G(n)^2-2G(n)G(n+1)+2G(n)G(n+2)
\qquad\qquad\label{eq:subword1}\\
&&\qquad\qquad+2(A-1)G(n+1)G(n+2)-(A-1)^2G(n+1)^2-G(n+2)^2\biggr) \nonumber
\end{eqnarray}
if $A>2$, and
\begin{eqnarray}
g(n)&=& n+2G(n)^2+7G(n)G(n+1)-5G(n)G(n+2)\label{eq:subword2}\\
&&\qquad\qquad\qquad\qquad\qquad-5G(n+1)G(n+2)+3G(n+1)^2+2G(n+2)^2\nonumber
\end{eqnarray}
if $A=2$.

In the case $A=2$, it is easy to show that $G(n)=F(n+3)-1$ for all $n$, where
$F(n)$ is the $n$th Fibonacci number. Consequently,
$G(n+2)-G(n+1)-G(n)=1$, and adding any multiple of the relation
\begin{equation}\label{eq:rel2}
(G(n+2)-G(n+1)-G(n)-1)^2=0
\end{equation}
to the right side of (\ref{eq:subword2}) gives another degree $2$ expression
of
type (\ref{eq:ans1}) or (\ref{eq:ans3}). Thus formula (\ref{eq:subword2}) is
not unique
within the class of formul\ae\ of type (\ref{eq:ans1}) or (\ref{eq:ans3}).
However, in
Section 5 we will show that, when $A>2$, formula (\ref{eq:subword1}) is
unique within this
class. When $A=2$, we show that all
relations are constant multiples of (\ref{eq:rel2}).

\subsection{Fibonacci power sums}
Theorems 1 and 2 imply that if the $F(j)$'s are the Fibonacci
numbers then for
each integer $p=1,2,\dots$ there is a formula
$$
f(n)=
\sum_{j=0}^{n-1}F(j)^p=\sum_{j=0}^p\Lambda_{p,j}F(n)^jF(n+1)^{p-j}+c_pn+d_p.
$$
Here is a brief table of values of these coefficients.

%\smallskip

\[\begin{array}{r|rrrrrrrrrr}
p&\Lambda_{p,0}&\Lambda_{p,1}&\Lambda_{p,2}&\Lambda_{p,3}&\Lambda_{p,4}&
\Lambda_{p,5}&\Lambda_{p,6}&\Lambda_{p,7}&c_p&d_p\\
\hline
&&&&&&&&&&\\
1&1&0&0&0&0&0&0&0&0&-1\\
&&&&&&&&&&\\
2&0&1&-1&0&0&0&0&0&0&0\\
&&&&&&&&&&\\
3&-\frac12&\frac32&0&-\frac32&0&0&0&0&0&\frac12\\
&&&&&&&&&&\\
4&0&\frac{2}{25}&-\frac{3}{25}&\frac{14}{25}&-\frac{19}{25}&0&0&0&\frac{6}{25}&0
\\
&&&&&&&&&&\\
5&\frac{7}{22}&-\frac{5}{22}&-\frac{15}{11}&\frac{10}{11}&\frac{15}{11}&
-\frac{15}{22}&0&0&0&-\frac{7}{22}\\
&&&&&&&&&&\\
6&0&\frac12&-\frac{5}{4}&0&\frac{5}{4}&
\frac{1}{2}&-1&0&0&0\\
&&&&&&&&&&\\
7&-\frac{139}{638}&\frac{763}{638}&-\frac{945}{638}&-\frac{350}{319}&\frac{105}{58}&
\frac{357}{319}&-\frac{105}{319}&-\frac{777}{638}&0&\frac{139}{638}
\end{array}\]

%\bigskip

\noindent
The resulting expressions for $f(n)$ turn out
to be unique within the class of type (\ref{eq:ans1}) or
(\ref{eq:ans3}) formul\ae\ when $p\not\equiv 0 \mbox{ mod } 4$.
When $p$ is a multiple of $4$ (for example, in the fourth line
of the table above) the
formul\ae\ are not unique, but are subject to
a one-parameter family of relations generated by powers of the degree-$4$
relation
$$
(F(n+1)^2 - F(n)^2 - F(n)F(n+1))^2 = 1\quad.
$$
We will establish these facts in Section 5.

\subsection{Generic power sums}

Consider the power sum
$$
f(n) = \sum_{k=0}^{n-1} F(k)^2,
$$
where $F(n)$
solves a linear recurrence
$$ F(n) = A F(n-1) + B F(n-2)
$$
with initial values $F(0)$ and $F(1)$, where $A$
and $B$ are sufficiently general to insure that, if $r_1$ and $r_2$ are
the associated roots, then $r_1$ and $r_2$ are distinct and
none of the monomials
$r_1^2,r_1 r_2$, and $r_2^2$ equals $1$. This is equivalent to assuming
simply that
$A^2+4B \not=0$, $A\not= \pm (B-1)$, and $B\not= -1$.

Using techniques introduced earlier, we can express $f(n)$ as a linear
combination of $F(n)^2, F(n)F(n+1), F(n+1)^2$, and $1$. The
solution may be computed explicitly in terms of $A,B,F(0)$, and
$F(1)$, and we find that $f(n)$ equals
\begin{equation}\label{generic}
\frac{(1-B-A^2(1+B)) F(n)^2 + (2AB )F(n)F(n+1) +
(1-B)F(n+1)^2 - K}{(A^2 - (B-1)^2)(B+1)}
\end{equation}
where $$K= (1-B-A^2(B+1))F(0)^2 + (2AB) F(0) F(1) + (1-B)F(1)^2.$$
In (\ref{generic}), we observe a curious phenomenon: since $F(n)$ depends on
$F(0)$ and $F(1)$, we might expect that our linear equations would have led
to a solution in which each of the
coefficients depends on $F(0)$ and $F(1)$. However, this
dependence appears only in the constant term. The next theorem demonstrates that such behavior is
typical for power sums of C-finite functions in which the terms in (\ref{eq:ans1}) containing $n$
are not present, i.e, in cases where no monomial in the roots equals $1$.
\begin{theorem}\label{inits}
Suppose that $\{F(n)\}_{n\geq 0}$ is a C-finite sequence determined by a recurrence of order $d$
together with initial conditions $F(0),F(1),\dots,F(d-1)$. Suppose that the recurrence
polynomial has distinct roots $r_1,\dots,r_d$, and suppose that no monomial of degree
$p$ in the $r_i$ equals $1$.
Let
$f(n) = \sum_{j=0}^{n-1} F(j)^p$, where $p$ is a positive integer, and let
\begin{equation}\label{initseq}
f(n) = \sum_{0\leq i_1,i_2,\dots ,i_d\leq d-1} \Lambda_{i_1,i_2,\dots ,i_d} F(n+i_1) F(n+i_2) \cdots F(n+i_d )\;+\;  K
\end{equation}
be the expansion of $f(n)$ obtained according to the method given
in Section 2. Then the
coefficients $\Lambda_{i_1,i_2,\dots ,i_d}$ in (\ref{initseq}) do not depend on $F(0),F(1),\dots,F(d-1)$.
\end{theorem}

\noindent
{\bf Proof: }
Suppose that $F(n) = \sum_{m=1}^d \lambda_m r_m^n$. Define
$$
{\bf X}(n) =
\left(\begin{array}{c}
\lambda_1 r_1^n \\
\lambda_2 r_2^n \\
\vdots \\
\lambda_d r_d^n
\end{array}\right)
\quad \mbox{and}\quad {\bf Y}(n) = \left(\begin{array}{c}
F(n) \\
F(n+1)\\
\vdots \\
F(n+d-1)
\end{array}\right).
$$
Then we have
\begin{equation}\label{invertv}
{\bf Y}(n) \;=\; {\bf V} \,{\bf X}(n)\quad\mbox{and}\quad {\bf X}(n) = {\bf V}^{-1} {\bf Y}(n)
\end{equation}
where $\bf V$ is a Vandermonde matrix in the $r_i$. It follows from
(\ref{invertv}) that the terms $\lambda_m r_m^n, 1 \leq m \leq d,$ can be expressed as linear combinations
of the functions $F(n+i)$, with coefficients that do not depend on
$F(0),F(1),\dots,F(d-1)$. Using the method of Section 2, we can compute
\begin{eqnarray}
f(n)=\sum_{j=0}^{n-1} F(j)^p &=&
\sum_{j=0}^{n-1} \bigg( \sum_{m=0}^d \lambda_m r_m^j \bigg)^p\nonumber\\
&=&
\sum_{j=0}^{n-1}\bigg( \sum_{0 \leq i_1,i_2,\dots,i_p\leq d}
\lambda_{i_1} \lambda_{i_2}\cdots  \lambda_{i_p}(r_{i_1} r_{i_2}\cdots  r_{i_p})^j\bigg)\nonumber\\
&=&
\sum_{0 \leq i_1,i_2,\dots,i_p\leq d} \lambda_{i_1} \lambda_{i_2}\cdots  \lambda_{i_p}
\frac{(r_{i_1} r_{i_2}\cdots  r_{i_p})^n-1}{(r_{i_1} r_{i_2}\cdots  r_{i_p})-1}\nonumber\\
&=&
\sum_{0 \leq i_1,i_2,\dots,i_p\leq d}
\frac{(\lambda_{i_1}r_{i_1}^n)(\lambda_{i_2}r_{i_2}^n)\cdots
(\lambda_{i_p}r_{i_p}^n)}{(r_{i_1}r_{i_2}\cdots r_{i_p}) -1} \;-\;K \label{final}
\end{eqnarray}
where $K$ is a constant. Using (\ref{invertv}), we can express all
of the terms in (\ref{final}) except $K$ as a linear combination
of monomials in the $F(n+i)$ with coefficients that do not depend
on $F(0),F(1),\dots,F(d-1)$, as claimed. $\Box$

We will return to this subject in Section 8, where we prove that a more general
version of Theorem \ref{inits} holds even when the roots $r_i$ are not distinct.

\subsection{Computational issues}

We have seen in the above theorems and corollaries that we can
decide the uniqueness of representations of certain sums in
closed form if we can decide whether or not the $N={p+d-1\choose
p}$ formally distinct monomials of degree $p$ in the roots
$r_1,\dots,r_d$ actually are all different, when evaluated as
complex numbers. We note here that there are various ways in which this can be done  without computing the roots.

For example, the elementary symmetric functions of these $N$
monomials in the $r_i$'s are symmetric functions in the $r_i$'s
themselves. Since any symmetric function of the roots of a
polynomial can be computed rationally in terms of its
coefficients, the same applies to these. Once the elementary
symmetric functions of the $N$ monomials of degree $p$ have been
computed, the discriminant of the polynomial whose coefficients
they are can be computed in the usual way. Thus, our condition on
the roots of $F$ can be tested without finding the roots. It would
be interesting to investigate in general this
``hyperdiscriminant'' of degree $p$ that is attached to a
polynomial $f$, particularly with regard to how it factors when
expressed in terms of the coefficients of $f$.

\section{Uniqueness and dimension: Fibonacci power sums}
In the next two sections, we investigate the uniqueness of the expansions
guaranteed by Theorems \ref{thm1} and \ref{thm2}.
Motivated by Example 3.3, we first consider this question for
expansions of the form
$$
\sum_{j=0}^p\Lambda_{p,j}F(n)^jF(n+1)^{p-j}
$$
and, more generally,
$$
\sum_{j=0}^p\Lambda_{p,j}F(n)^jF(n+1)^{p-j}+c_pn+d_p,
$$
where $F(n)$ denotes the $n$th Fibonacci
number. In Section 5 we develop tools to help answer
these questions for
more general linear recurrences, and for other summations such
as those arising in Examples 3.1 and 3.2. The techniques in
these two sections can be viewed as refinements and extensions
of the ideas introduced in Section 2 to prove
Theorems \ref{thm1}, \ref{thm2},
and \ref{lem3}.

\begin{theorem}\label{th1}
Let $\bf V = \bVinf$ denote the vector space of complex-valued functions on
$\{0,1,2,\dots\}$, and let $\bf W_p$ denote the subspace of $\bf V$ spanned by
functions of the form $F(n)^i F(n+1)^{p-i}$ for
$i=0,\dots,p$, and let $\bf W_p^{++}$ denote the subspace spanned by
the same monomial expressions together with with the functions $g(n) =n$ and
$h(n)=1$. Then
\begin{enumerate}
\item
$\dim({\bf W_p}) = p+1$, and
\item
$\dim({\bf W_p^{++}}) =
\left\{
\begin{array}{cl}
p+2, & \mbox{ if p is divisible by $4$;}\\
p+3, & \mbox{ otherwise.}
\end{array}\right.
$
\end{enumerate}
\end{theorem}

\begin{corollary} The functions
$F(n)^i F(n+1)^{p-i} , 1 \leq i \leq p$ are linearly independent, and the set
$$\{F(n)^i F(n+1)^{p-i} \}_{1 \leq i \leq p} \cup \{n,1\}$$ is linearly
independent unless $p$ is divisible by $4$, in which case there
is a single relation among its elements.
\end{corollary}
\noindent{\bf Proof:} Let $r_1= (1+\sqrt{5})/2$ and $r_2 =
(1-\sqrt{5})/2$ denote the roots of the Fibonacci recurrence
polynomial. As noted earlier in the proof of Theorem 1, $r_1^n$ and $r_2^n$ may
be expressed
as linear combinations of $F(n)$ and $F(n+1)$ and vice versa.
Consequently, $\bf W_p$ is the linear span of
$r_1^{ni} r_2^{n(p-i)}, i=0,\dots, p$, and to prove statement (a) it
suffices
to show that these functions are linearly independent. But this follows
immediately
from the fact that the numbers $r_1^{i} r_2^{(p-i)}$
are distinct, for $i=0,\dots,p$.

To prove part (b), consider the $(p+3) \times (p+3)$ matrix $M_p$ whose
$i$th column is equal to the vector
$(1,\theta_i,\theta_i^2,\dots,\theta_i^{p+2})$, where
$\theta_i=r_1^{i} r_2^{(p-i)}, i=0,\dots,p$, and whose last two columns are
the vectors $(1,1,\dots,1)$ and $(0,1,\dots,p+2)$. For example, when $p=2$ we
have
$$
M_2 \;=\;
\left(
\begin{array}{ccccc}
1 & 1 & 1 & 1 & 0\\
r_2^2 & r_1 r_2 & r_1^2 & 1 & 1\\
r_2^4 & r_1^2 r_2^2 & r_1^4 & 1 & 2\\
r_2^6 & r_1^3 r_2^3 & r_1^6 & 1 & 3\\
r_2^8 & r_1^4 r_2^4 & r_1^8 & 1 & 4
\end{array}
\right)
$$
Note that $\det M_p$ is the derivative at $t=1$ of the $(p+3) \times (p+3)$
Vandermonde determinant
$\det M_p(t)$, where $M_p(t)$ is the matrix whose
first $p+2$ columns are the same as those of $M_p$, and whose last column is
$(1,t,t^2,\dots,t^{p+2})$. We have
\begin{eqnarray*}
\det M_p &=& \frac{d}{dt}\; \det M_p(t) \;\bigg|_{t=1} \\
&=& \frac{d}{dt}\; \biggl( \prod_{0\leq i<j\leq p} (\theta_j-\theta_i)
\prod_{0\leq i \leq p} (1-\theta_i)  \prod_{0\leq i \leq p}
(t-\theta_i)\;(t-1)\biggr) \;\bigg|_{t=1}\\
&=& \prod_{0\le i<j\le p} (\theta_j-\theta_i) \prod_{0\le i \le p}(1-\theta_i)^2
\nonumber
\end{eqnarray*}
It follows that $\det M_p=0$ only when $t=1$ is a multiple root of $\det
M_p(t)$,
i.e., $r_1^i r_2^{p-i}= 1$ for some $i$.  Using the fact that $r_1 r_2
= -1$, it is easy to
show that this property holds if and only if $p$ is a multiple of $4$. Thus,
when $p$ is not a multiple of $4$,
the columns of $M_p$ are linearly independent and we have
$\dim ({\bf W_p^{++}})= p+3$.

If $p$ is a multiple of $4$, then $M_p$ contains exactly two columns of
$1$s. If one of these
columns is
suppressed, the argument just given shows that the remaining columns are
linearly independent.
Hence $\mbox{rank}(M_p)= p+2$ and $\dim ({\bf W_p^{++}})\geq p+2$. Since the
dimension
is clearly at most $p+2$ in this case, the theorem is proved. $\Box$

\section{Uniqueness and dimension: other recurrences with distinct roots}

Analogs of Theorem \ref{th1} hold for more general recurrences with distinct roots,
but the exact statements depend on properties of
the associated roots. The following theorem
concerns relations among monomials of type (\ref{eq:ans3}),
and allows precise dimension computations in many cases.

\begin{theorem}\label{thm:dimw}
Let $F(n)$ be a solution to a linear recurrence of order $d$ whose
associated roots
$r_1,r_2,\dots,r_d$
are distinct, and let $p$ and $q$ be distinct
positive integers. Let ${\bf W_p}$ denote the subspace of ${\bf V}=\bVinf$ spanned
by all degree $p$ monomials of the form
$$
F(n)^{i_1} F(n+1)^{i_2}\cdots F(n+d-1)^{i_d}
$$
where $i_1+i_2 + \cdots i_d=p$ and $i_j \geq 0$ for all $j$. Let
${\bf W_q^+}$ denote the subspace spanned by all degree $q$
monomials
\begin{eqnarray*}
F(n)^{i_1} F(n+1)^{i_2}\cdots F(n+d-1)^{i_d} &\mbox{and}&\\
n F(n)^{i_1} F(n+1)^{i_2}\cdots F(n+d-1)^{i_d}
\end{eqnarray*}
where $i_1+i_2 + \cdots i_d=q$ and $i_j \geq 0$ for all $j$.
And, finally, let ${\bf W_{p,q}^{++}}= {\bf W_p}+{\bf W_q^+}$
denote the subspace spanned by all of the above monomials.
Then
$$\dim({\bf W_p})= |S_p|\;,
\quad \dim({\bf W_q^+})= 2|S_q|\;,
\quad \mbox{and}\quad \dim({\bf W_{p,q}^{++}})= |S_p|+ 2 |S_q| - |S_p \cap
S_q|
$$
where
$S_p=\{r_1^{i_1}r_2^{i_2}\cdots r_d^{i_d}\;|\; i_1+i_2 + \cdots i_d=p\} $
and
$S_q=\{r_1^{i_1}r_2^{i_2}\cdots r_d^{i_d}\;|\; i_1+i_2 + \cdots i_d=q\} $
are the sets of
monomials in the $r_i$ of degrees $p$ and $q$, respectively, both
viewed as subsets of the complex numbers.
\end{theorem}

\begin{corollary}\label{corw}
The sets of monomials generating ${\bf W_p}$, ${\bf W_q^+}$, and ${\bf
W_{p,q}^{++}}$,
respectively, are linearly independent if and only if evaluations of
formally distinct monomials
in the sets
$S_p$, $S_q$ and $S_p \cup S_q$ yield distinct complex numbers.
\end{corollary}

\noindent
The proof of Theorem \ref{thm:dimw}
is analogous to that given for Theorem \ref{th1}, but more
careful analysis is required. First consider the case of ${\bf W_p}$. As noted in
Section 2, each of the functions
$F(n), F(n+1),\dots,F(n+d-1)$ lies in the linear span of
$r_1^n,r_2^n,\dots,r_d^n$, and
conversely. Hence ${\bf W_p}$ is spanned by the set
$\{\theta_1^n,\theta_2^n,\dots,\theta_{m(p,d)}^n\}$, where
$m(p,d)= {{p+d-1}\choose{p}}$ and
the $\theta_j$ range over the $m(p,d)$ formally distinct monomials of degree
$p$
in $r_1,r_2,\dots,r_d$.
Similar reasoning shows that ${\bf W_q^+}$ is spanned by the $2m(q,d)$
functions
\begin{eqnarray}
&\psi_1^n,\psi_2^n,\dots,\psi_{m(q,d)}^n&\label{eq:wq}\\
&n\psi_1^n,n\psi_2^n,\dots,n\psi_{m(q,d)}^n&,\nonumber
\end{eqnarray}
where
the $\psi_j$ range over all formally distinct monomials of degree
$q$ in $r_1,\dots,r_d$, and finally,
${\bf W_{p,q}^{++}}$ is spanned by the $m(p,d)+2m(q,d)$ functions
\begin{eqnarray}
&\theta_1^n,\theta_2^n,\dots,\theta_{m(p,d)}^n&\nonumber\\
&\psi_1^n,\psi_2^n,\dots,\psi_{m(q,d)}^n&\label{eq:wpq}\\
&n\psi_1^n,n\psi_2^n,\dots,n\psi_{m(q,d)}^n,& \nonumber
\end{eqnarray}
where $\theta_i$ and $\psi_j$ are defined as above.
Theorem \ref{thm:dimw} is now an immediate consequence of the following Lemma.
\begin{lemma}\label{standard}
Let $\omega_1,\omega_2,\dots,\omega_m$ be complex numbers, and let $a_1,a_2,\dots,a_m$ be positive integers.
Then the functions
$$
n^j \omega_i^n \qquad 1 \le i \le m,\; 0 \le j \le a_i -1
$$
are linearly independent if and only if the $\omega_i$ are distinct.
\end{lemma}
Lemma \ref{standard} is a standard component of the classical theory of finite
difference equations (e.g., \cite{jordan}, Chapter 11), indeed it is the justification for the usual method of
solution of such equations. It is easy to give a direct proof via
generating functions, or, alternatively, one can give a
Vandermonde-type proof based on the following elegant
determinant formula (\cite{fh}, but also see \cite{k} for an extensive
history of this formula).

\begin{theorem}\label{nicematrix}
Let $x_1,x_2,\dots,x_n$ be indeterminates, and let $a_1,a_2,\dots,a_n$ be
positive integers with $\sum_i a_i = N$.
For all $t$, and for any integer $k\geq 1$, let
$$\rho_N(t,k) = \frac{d^{k}}{dt^{k}}(1,t,t^2,\dots,t^{N-1})$$
Let $M(a_1,a_2,\dots,a_n)$ be the $N \times N$ matrix whose first $a_1$ rows
are
$\rho_N(x_1,0),\dots,\rho_N(x_1,a_1-1)$, and whose next $a_2$ rows are
$\rho_N(x_2,0),\dots,\rho_N(x_2,a_2-1)$, and so forth.
Then
$$
\det M(a_1,\dots,a_n)\;=\;
\prod_{i=1}^n (a_i-1)!!!
\prod_{1\leq i < j \leq n} (x_j-x_i)^{a_i a_j}
$$
where $k!!!$ denotes $1! 2! \cdots k!$ and $0!!!= 1$.
\end{theorem}
For example,
$$
M(1,2,3) =
\left(
\begin{array}{cccccc}
1 & x_1 & x_1^2 & x_1^3 & x_1^4 & x_1^5 \\
1 & x_2 & x_2^2 & x_2^3 & x_2^4 & x_2^5 \\
0 & 1 & 2x_2 & 3x_2^2 & 4x_2^3 & 5x_2^4\\
1 & x_3 & x_3^2 & x_3^3 & x_3^4 & x_3^5 \\
0 & 1 & 2x_3 & 3x_3^2 & 4x_3^3 & 5x_3^4\\
0 & 0 & 2 & 6x_3 & 12x_3^2 & 20x_3^3
\end{array}
\right)
$$
and
$$
\det M(1,2,3) \;=\; 2 (x_2-x_1)^2 (x_3-x_1)^3 (x_3-x_2)^6.
$$

Theorem \ref{thm:dimw} describes relations among closed form expressions of
type (\ref{eq:ans3}), but the proof also yields similar results
for expressions of type (\ref{eq:ans1}).
\begin{corollary}\label{cortype3}
Let $F(n)$ be a solution to a linear recurrence of order $d$ whose associated
roots are distinct. Let ${\bf W^*_d}$ denote the space spanned
by monomial functions of type (\ref{eq:ans1}). Then $\dim {\bf W^*_d} =
|S|+2|T|$, where $S$
is the set of all monomials of the form
$t_1^{a_1} t_2^{a_2} \cdots t_k^{a_k}$ and $T$ is the set of
monomials of the form
$t_1^{a_1+b_1} t_2^{a_2+b_2} \cdots t_k^{a_k+b_k}$
and, for each $i$, $t_i$ is one of the roots $r_1,r_2,\dots,r_d$.
The monomial functions of type
(\ref{eq:ans1}) are linearly
independent if and only if formally distinct monomials in $S\cup T$ correspond
to distinct complex numbers.
\end{corollary}
We omit the proof, which is analogous to that of the proof of Theorem
\ref{thm:dimw}.
We note that the set
$S\cup T$ in Corollary \ref{cortype3} is a subset of the set
$S_p\cup S_q$ appearing in Corollary \ref{corw},
and thus we obtain the following result.
\begin{corollary}\label{cor13}
Under the assumptions of Corollary \ref{cortype3}, if the monomial functions of type
(\ref{eq:ans3}) are
linearly independent, then so are the monomial functions of type
(\ref{eq:ans1}).
\end{corollary}

We will now apply these results to some of the formul\ae\ in Sections 3.1
and 3.2.

\begin{corollary}
For the Fibonacci sum $f(n)$ appearing in (\ref{eq:coll}), equation
(\ref{eq:coll2})
gives the unique $F$-closed formula of type (\ref{eq:ans1}) and
(\ref{eq:coll2a}) gives the unique $F$-closed
formula of type (\ref{eq:ans3}).
\end{corollary}

\noindent
{\bf Proof:} By Theorem \ref{thm:dimw} and Corollary \ref{cor13}, we need only
check that, if $r_1$ and
$r_2$ denote the roots of the Fibonacci recurrence, then
$$
r_1^2, r_1 r_2, r_2^2, r_1^3, r_1^2 r_2, r_1 r_2^2, \mbox{ and } r_2^3
$$
are distinct real numbers. This is an elementary calculation.

\begin{corollary}
For the sum $g(n)=\sum_{j=0}^{n-1}G(j)^2$ arising in the subword
avoidance problem with $A=2$, solutions $g(n)$ of type
(\ref{eq:ans3})
are all given by
(\ref{eq:subword2}) plus constant multiples of relation (\ref{eq:rel2}).
\end{corollary}

\noindent
{\bf Proof: } The roots of the recurrence equation $t^3 -2 t^2 + 1=0$ are
$r_1, r_2, r_3$, where $r_1$ and $r_2$ are roots of the Fibonacci recurrence
and $r_3=1$. By Theorem \ref{thm:dimw}, the dimension of the space
${\bf W_{2,0}^{++}}$ spanned by the six degree-$2$
monomials in $G(n), G(n+1)$ and $G(n+2)$ together with $1$ and $n$ is equal to
$|S_2|+ 2|S_0|- |S_2\cap S_0|$, where
$S_2=\{r_1^2,r_2^2,1,r_1r_2,r_1,r_2\}$,
and $S_0=\{1\}$. Elementary calculation shows that this dimension is equal
to $7$, hence the monomials generating ${\bf W_{2,0}^{++}}$ are linearly
independent apart from a one-parameter family of relations.

Next we consider the case $A=3$, as a warmup for the general case $A>2$.

\begin{corollary}
For the sum $g(n)=\sum_{j=0}^{n-1}G(j)^2$ arising in the subword
avoidance problem with $A=3$, formula
(\ref{eq:subword1})
gives the unique $G$-closed formula of type (\ref{eq:ans1}).
\end{corollary}

\noindent
{\bf Proof: } Here the recurrence equation is $t^3 - 3t^2 +1=0$, which has
roots
$r_1=1 + \eta + \eta^{17}, r_2 = 1 + \eta^7 + \eta^{11}, r_3 = 1 + \eta^5 +
\eta^{13}$,
where $\eta = e^{2\pi i/18}$ is an $18$th root of unity. Again, by Theorem
\ref{thm:dimw}, the
dimension of the space of monomials ${\bf W_{2,0}^{++}}$ is equal to
$|S_2|+ 2|S_0|- |S_2\cap S_0|$, where
$S_2=\{r_1^2,r_2^2,r_3^2,r_1r_2,r_1r_3,r_2r_3\}$,
and $S_0=\{1\}$. A slightly less elementary calculation shows that
the formal monomials in $S_2\cup S_0$ are distinct, so that
$\dim({\bf W_{2,0}^{++}})=8$ and the monomial functions generating
${\bf W_{2,0}^{++}}$ are linearly independent.

\begin{corollary}\label{cor:niceproof}
For the more general power sum $g(n)=\sum_{j=0}^{n-1}G(j)^p$ arising in the
subword
avoidance problem, with $p>0$ and any $A>2$, solutions
of type (\ref{eq:ans1})
are unique if and only if $p \not\equiv 0 \mbox{ mod } 6$.
\end{corollary}

\noindent
{\bf Proof: } An argument analogous to the calculation in Section 3.2 shows
that formul\ae\
of type (\ref{eq:ans1}) exist expressing $g(n)$ as linear combinations of
monomials in
$G(n), G(n+1)$, and $G(n+2)$ of degree $p$, together with $1$ and $n$. We
need to compute
the dimension of ${\bf W_{p,0}^{++}}$, which by Theorem \ref{thm:dimw} is
equal to
$|S_p|+ 2|S_0|- |S_p\cap S_0|$, where $S_0=\{1\}$ and $S_p$ is the set of
all degree-$p$
monomials in $r_1, r_2$, and $r_3$, where $r_1, r_2$, and $r_3$ are roots of
the recurrence equation $t^3-At^2 + 1= 0$.

If $p$ is not divisible by $6$, the proof will be complete if we can
show that formally distinct monomials in $S_p$
evaluate to distinct complex (actually real) numbers, and none of them
equals $1$.
Suppose that $r_1^{e_1}r_2^{e_2}r_3^{e_3}=r_1^{f_1}r_2^{f_2}r_3^{f_3}$,
where
$\sum e_i = \sum f_i = p$ and $e_i\not=f_i$ for some $i$. Then by
cancellation we
obtain the relation $r_i^{u_i}=r_j^{u_j}r_k^{u_k}$ for some rearrangement of
the indices,
with $u_i,u_j,u_k\geq 0$ and at least one of these exponents positive.
Using the relation $r_1r_2r_3 =-1$, if necessary, to eliminate one of the
roots,  we obtain (after possibly reindexing), $r_i^{v_i} = \pm r_j^{v_j}$
with
$v_i,v_j\geq 0$ and at least one of these exponents positive.

It is a straightforward exercise to show that the roots
$r_1,r_2$ and $r_3$ are all real, and that,
if they are arranged in decreasing order, then
$r_1>1, 0<r_2<1$, and $-1 < r_3 < 0$. From elementary Galois theory we know that there exists an automorphism
$\Phi$
of the field $K ={\bf Q}(r_1,r_2,r_3)$ such that
$\Phi: r_1 \mapsto r_2 \mapsto r_3 \mapsto r_1$, i.e., it permutes the roots
cyclically. Hence the equation $r_i^{v_i} = \pm r_j^{v_j}$
holds for all three cyclic permutations of the roots. At least one of these
equations leads to a contradiction, since $|r_1|>1$ and $|r_2|,|r_3|<1$.
This proves
that formally distinct monomials are distinct, and it remains to show that
none can equal $1$.

Suppose that $r_1^{e_1}r_2^{e_2}r_3^{e_3}=1$, and the exponents $e_i$ are
not all
equal. Applying the identity $r_1r_2 r_3 = -1$ we obtain a relation of the
form $r_i^{u_i} r_j^{u_j} = \pm 1$ for some pair of distinct $i,j$, with
$u_i,u_j\geq 0$ and at least one positive. Again, this relation holds for
all cyclic
permutations of the indices, and consideration of absolute values leads to a
contradiction in at least one case. Consequently, we must have
$e_1=e_2=e_3=e$ for some $e$. From the relations $r_1 r_2 r_3=-1$ and
$(r_1r_2r_3)^e=1$
we conclude that $e$ is even,
which implies that $p$ is a multiple of $6$. This completes the proof
that monomials in the $G$ are linearly independent when $p \not\equiv 0
\mbox{ mod } 6$.
When $p=6m$ the relation $(r_1 r_2 r_3)^{2m}=1$ gives relations in the $G$
of
degree $6$, and so the proof of Corollary \ref{cor:niceproof} is complete.
$\Box$

\section{The general case of multiple roots}
In this section we show the result of dropping the assumption of distinct roots.
We also consider a somewhat more general summation problem, viz.
\begin{equation}
\label{eq:gensum}
f(n)=\sum_{j=0}^{n-1} F_1(a_1n+b_1j+c_1)F_2(a_2n+b_2j+c_2)\cdots
F_k(a_kn+b_kj+c_k),
\end{equation}
in which the factors of the summand may be different $C$-finite functions.
The analysis in this general case is similar to that in the case of distinct roots,
but some additional machinery is required. The main result is the following,
which is a generalization and also a strengthening of Theorem \ref{thm1}.

\begin{theorem}
\label{th:mpl} Let $F_1,F_2,\dots,F_k$ be given $C$-finite
sequences. Suppose that, for each $i=1,\dots,k$, $F_i(n)$
satisfies a recurrence of minimal degree $D(i)$ whose polynomial equation
has $d(i)$
distinct roots. Denote these roots by $r_1^{(i)},\dots,r_{d(i)}^{(i)}$, and let
$e^{(i)}_1,\dots,e^{(i)}_{d(i)}$ be their respective
multiplicities, so that $D(i)=\sum_j e^{(i)}_j$. Finally, let
$\Delta(i)=\max_{1\le j\le d(i)}(e^{(i)}_j-1)$. Then the sum
$f(n)$, of (\ref{eq:gensum}), can be expressed as a linear
combination of the monomials
\begin{equation}\label{eq:form1}
F_1((a_1+b_1)n+i_1)\dots F_k((a_k+b_k)n+i_k)\quad(0\le i_{\nu}\le D(\nu)-1;1\le \nu\le k),\ \mathrm{and}
\end{equation}
\begin{equation}
\label{eq:form2}
\psi_{i_1,\dots,i_k}(n)F_1(a_1n+i_1)\dots F_k(a_kn+i_k)\qquad\quad\qquad(0\le i_{\nu}\le D(\nu)-1;1\le \nu\le k),
\end{equation}
in which $\psi_{i_1,\dots,i_k}(n)$ is a polynomial whose degree is bounded above by
\begin{itemize}
\item zero, i.e., the factor $\psi$ can be omitted, if for all
sequences $(m_1,m_2,\dots,m_k)$ with $1 \leq m_i \leq d(i)$,
we have $\prod_{i=1}^k(r_{m_i}^{(i)})^{b_i}\neq 1$, and
\item $1+\sum_i \{\Delta(i): a_i=0\}$ otherwise.
\end{itemize}
\end{theorem}

\noindent{\bf Proof: }We have, for the sum $f(n)$ of (\ref{eq:gensum}),
\begin{equation}
\label{eq:expd}
f(n)=\sum_{j=0}^{n-1}\prod_{\ell=1}^k\sum_{m=1}^{d(\ell)}\sum_{h=0}^{e_m^{(\ell)}-1}
\lambda_{m,h}^{(\ell)}(a_{\ell}n+b_{\ell}j+c_{\ell})^{h}(r_m^{(\ell)})^{a_{\ell}n+b_{\ell}j+c_{\ell}},
\end{equation}
where the coefficients $\lambda_{m,h}^{(\ell)}$ are defined by the form of the $F_i$'s, namely
\[F_i(n)=\sum_{m=1}^{d(i)}\sum_{h=0}^{e_m^{(i)}-1}\lambda_{m,h}^{(i)}n^{h}(r_m^{(i)})^n.\]
If we expand the product and the two inner sums
in (\ref{eq:expd}), we find that a typical term is of the form
\[K(a_1n+b_1j+c_1)^{h_{1}}(r_{m_1}^{(1)})^{a_1n+b_1j+c_1}\dots (a_kn+b_kj+c_k)^{h_{k}}(r_{m_k}^{(k)})^{a_kn+b_kj+c_k},\]
where $0\le h_i \le e_{m_i}^{(i)}-1$ for $1\le i\le k$.
If we write
$$(a_i n + b_i j + c_i)^{h_i} = (a_i(n-j) + (a_i+b_i)j + c_i)^{h_i}
$$
and further expand each of these factors,
we find that our typical term now can be expressed as
\begin{equation}
\label{eq:ofform}
Kj^q(n-j)^r\left((r_{m_1}^{(1)})^{b_1}\dots
(r_{m_k}^{(k)})^{b_k}\right)^j\left((r_{m_1}^{(1)})^{a_1}\dots (r_{m_k}^{(k)})^{a_k}\right)^n,
\end{equation}
in which
\begin{equation}
\label{eq:bounds}
q\le\sum\{h_i:a_i+b_i\neq 0\},\quad r\le \sum \{h_i:a_i\neq 0\},\mbox{    and      }
q+r \le \sum_i h_i \le \sum_i (e_{m_i}^{(i)}-1).
\end{equation}
At this point we need the following result.
\begin{lemma} We have
\[
\sum_{j=0}^{n-1}j^a(n-j)^bx^j=
\begin{cases}
{P_a(n)x^n+P_b(n)},&\text{if $x\neq 1$;}\\
{P_{a+b+1}(n)},&\text{if $x=1$,}
\end{cases}
\]
where $P_s(n)$ denotes a generic polynomial of degree $s$, whose coefficients may depend on $x$.
\end{lemma}
\noindent{\bf Proof.} Suppose that $x\not=1$. If $b=0$, we have
\[
\sum_{j=0}^{n-1}j^ax^j=\left(x\frac{d}{dx}\right)^a\sum_{j=0}^{n-1}x^j
=\left(x\frac{d}{dx}\right)^a\left(\frac{x^n-1}{x-1}\right),
\]
which is of the form stated when $b=0$. For $b>0$ we have
\begin{eqnarray*}
\sum_{j=0}^{n-1}j^a(n-j)^bx^j&=&\left(x\frac{d}{dx}\right)^a
\sum_{j=0}^{n-1}(n-j)^bx^j=\left(x\frac{d}{dx}\right)^a\sum_{j=1}^{n}j^bx^{n-j}=
\left(x\frac{d}{dx}\right)^ax^n\sum_{j=1}^{n}j^bx^{-j}\\
&=&\left(x\frac{d}{dx}\right)^ax^n\left((P_b(n)+n^b)x^{-n}+K\right)=
\left(x\frac{d}{dx}\right)^a\left(P_b(n)+n^b+Kx^n\right)
\end{eqnarray*}
which is evidently of the desired form. The case $x=1$ is elementary. $\Box$

If we sum the typical term (\ref{eq:ofform}) over $j=0,..,n-1$ and use the Lemma, we find that
the overall sum $f(n)$ is
a sum of expressions of the form
\begin{eqnarray}
\label{eq:nxform}
&&\sum_{j=0}^{n-1}Kj^q(n-j)^r\left((r_{m_1}^{(1)})^{b_1}\dots
(r_{m_k}^{(k)})^{b_k}\right)^j\left((r_{m_1}^{(1)})^{a_1}\dots
(r_{m_k}^{(k)})^{a_k}\right)^n\nonumber\\
&=&K\left((r_{m_1}^{(1)})^{a_1}\dots (r_{m_k}^{(k)})^{a_k}\right)^n\times
\begin{cases}
P_q(n)\left((r_{m_1}^{(1)})^{b_1}\dots (r_{m_k}^{(k)})^{b_k}\right)^n+P_r(n),& \text{if}\ \Theta\neq 1;\\
P_{q+r+1}(n), & \text{if}\ \Theta=1.
\end{cases}\nonumber\\
&=&
\begin{cases}
KP_q(n)\left((r_{m_1}^{(1)})^{a_1+b_1}\dots (r_{m_k}^{(k)})^{a_k+b_k}\right)^n+KP_r(n)
\left((r_{m_1}^{(1)})^{a_1}\dots (r_{m_k}^{(k)})^{a_k}\right)^n,& \text{if}\ \Theta\neq 1;\\[+6pt]
KP_{q+r+1}(n)\left((r_{m_1}^{(1)})^{a_1}\dots (r_{m_k}^{(k)})^{a_k}\right)^n, & \text{if}\ \Theta=1,
\end{cases}
\end{eqnarray}
where $\Theta=(r_{m_1}^{(1)})^{b_1}\dots (r_{m_k}^{(k)})^{b_k}$,
$q$ and $r$ satisfy the bounds given in (\ref{eq:bounds}), and $P_s(n)$ denotes a generic polynomial of degree $s$.

Considering each of the three terms appearing in the last member of (\ref{eq:nxform}), we first have
\begin{eqnarray*}
K P_q(n)\left((r_{m_1}^{(1)})^{a_1+b_1}\dots (r_{m_k}^{(k)})^{a_k+b_k}\right)^n&=&KP_{q}(n)
\prod_{\{i:a_i+b_i\neq 0\}}(r_{m_i}^{(i)})^{(a_i+b_i)n}\\
&=&\left(\sum_{j=0}^{q}\beta_jn^j\right)\prod_{\{i:a_i+b_i\neq 0\}}(r_{m_i}^{(i)})^{(a_i+b_i)n},
\end{eqnarray*}
say.
Since $q\le \sum_i\{h_i:a_i+b_i\not=0\} \le \sum_i\{e_{m_i}^{(i)}-1: a_i+b_i \not= 0\}$,
each exponent $j$ in the range $0\le j\le q$ can be
written (in many ways) as $j=j_1+j_2+\dots+j_k$, where
$0\le j_i \le e_{m_i}^{(i)}-1$ and $j_i = 0$ if $a_i+b_i=0$.
Hence the last member above may be expressed as
\begin{equation}
\label{eq:nform}
\sum_{j=0}^{q}\beta_j\prod_{\{i:a_i+b_i\neq 0\}}n^{j_i}(r_{m_i}^{(i)})^{n(a_i+b_i)}
\qquad(0\le j_i\le e_{m_i}^{(i)}).
\end{equation}
Now we observe that the solution space of the recurrence
satisfied by $F_i(n)$ has dimension $D(i)$, and that the $D(i)$ shifted
sequences $\{F_i(n)\},\{F_i(n+1)\},\dots,\{F_i(n+D(i)-1)\}$ are linearly independent,
since in the contrary case the function $F_i$ would satisfy a recurrence of degree $<D(i)$.
Consequently these $D(i)$ sequences are a basis for the solution space, and therefore each of the
functions
$n^j(r_{m}^{(i)})^n, j=0,1,\dots,e_m^{(i)}-1,m=1,\dots,d(i)$
can be written as a linear combination of $F_i(n),F_i(n+1),\dots,F_i(n+D(i)-1)$.

Thus we return to our general term (\ref{eq:nform}) and we
replace each of the monomials of the form $n^jr_{m}^{n(a+b)}$ by
such a linear combination of functions of the form
$F((a+b)n),F((a+b)n+1),\dots,F((a+b)n+D(i)-1)$, and expand everything
again. Now our general term is of the form (\ref{eq:form1}) in
the statement of the theorem, which concludes the
treatment of the first term in the final display of equation (\ref{eq:nxform}).
The second term,
$$
KP_r(n)
\left((r_{m_1}^{(1)})^{a_1}\dots (r_{m_k}^{(k)})^{a_k}\right)^n,
$$
may be handled similarly. In the third case, which occurs when $\Theta = 1$, powers of $n$
in $P_{q+r+1}(n)$ can be redistributed as in (\ref{eq:nform}) provided that each corresponding
$a_i$ is nonzero. If all $a_i \ne 0$, a
term of degree at most $1$ remains, since $q+r \le \sum_i e_{m_i}^{(i)}-1$; more generally, the residual exponent is
at most equal to $1 + \sum_i \{ \Delta(i): a_i = 0\}$.
$\Box$

By keeping track of the number of terms of the form $n^j \prod
r_m^{an}$ and $n^j \prod r_m^{(a+b)n}$ being summed in
(\ref{eq:nxform}) we can give a bound on the C-finite degree of
$f(n)$, generalizing Theorem \ref{lem3} to the case of multiple roots
and different $F_i$.
\begin{corollary}\label{gendegbound}
Suppose that $F_1, F_2,\dots, F_k$ are as defined in Theorem \ref{th:mpl}.
Let $M$ denote the number of algebraically
distinct monomials generated by (\ref{eq:form1}) and
(\ref{eq:form2}), as defined as in the paragraph preceding Theorem
\ref{lem3}. Then the sum $f(n)$ in (\ref{eq:gensum}) is C-finite,
of degree at most $M$.
The coefficients of the monomials expressing $f(n)$ as a linear combination
of those monomials
can be found by solving equations involving at most $M$ values of $f(n)$. If a solution
is valid for the first
$M$ values of $n$, then it is valid for all values of $n$.
\end{corollary}
\noindent
{\bf Proof:} Assume first that no
expression of the form $\Theta = \prod_{\ell=1}^k (r_{m_\ell}^{(\ell)})^{b_\ell}$ equals $1$.
Define $Q=\{\ell\;|\;a_\ell+b_\ell\not=0\}$ and
$R=\{\ell\;|\;a_\ell\not=0\}$.
The proof of Theorem \ref{th:mpl} shows that $f(n)$ is in the linear span of the set
$U_Q\cup U_R$, where $U_Q$ is the set of functions of the form
\begin{equation}\label{abprod}
\prod_{\ell\in Q} n^{j_\ell}(r_{m_\ell}^{(\ell)})^{(a_\ell+b_\ell)n}
\end{equation}
where $1 \le m_\ell \le d(\ell)$ and $0 \le j_\ell \le e_{m_\ell}^{(\ell)}-1$, and $U_R$ is
the set of functions of the form
\begin{equation}\label{aprod}
 \prod_{\ell\in R} n^{j_\ell}(r_{m_\ell}^{(\ell)})^{a_\ell n}
\end{equation}
where $1 \le m_\ell \le d(\ell)$ and $0 \le j_\ell \le  e_{m_\ell}^{(\ell)}-1$.
We claim that $|U_Q\cup U_R| \le M$. The argument is similar to
the one used in Section 2 to prove Theorem \ref{lem3} in the distinct root case: if $\theta \ne 0$, we define a map
$F_\ell(\theta n + i) \mapsto n^j (r_{m}^{(\ell)})^{\theta n}$, where $i \mapsto
(j,m)$ is some enumeration of the
$D(\ell)$ pairs with $1 \le m \le d(\ell), 0 \le j \le e_{m_\ell}^{(\ell)}-1$.
This map extends to a well-defined surjective map from the set of
equivalence classes of monomials
of type (\ref{eq:form1}) and (\ref{eq:form2}) to the set $U_Q \cup U_R$. The maximum C-finite degree of
any function in $\langle U_Q \cup U_R \rangle$ is equal to the dimension of that space, and
hence the C-finite degree of $f(n)$ is less than or equal to $M$.

Next suppose that there exist expressions of the form
$\Theta = \prod_{\ell=1}^k (r_{m_\ell}^{(\ell)})^{b_\ell}$ equal to $1$.
The proof of Theorem \ref{th:mpl} shows that $f(n)$ is in the linear span of the set
$U_Q\cup W_R$, where $U_Q$ is as defined above, and
$W_R$ is
the set of functions of the form
\begin{equation}\label{rhterm}
n^h \prod_{\ell\in R}(r_{m_\ell}^{(\ell)})^{a_\ell n}
\end{equation}
where $1 \le m_\ell \le d(\ell)$ and
$0 \le h \le H$, where
\begin{equation}\label{h0def}
H= 1 + \sum_{\ell \in  R} (e_{m_\ell}^{(\ell)}-1)+ \sum_{\ell \not\in R} (e_{m^*_\ell}^{(\ell)}-1)
\end{equation}
and, for a given sequence $m_\ell$ with $\ell\in R$, the $m^*_\ell$ are
chosen to maximize $H$ over all sequences $(m_1,\dots,m_k)$ such that
$\prod_{\ell\in R} (r_{m_\ell}^{(\ell)})^{b_\ell} \prod_{\ell\not\in R} (r_{m^*_\ell}^{(\ell)})^{b_\ell}= 1$.

As in the first case,
the map $F_\ell (\theta n + i)\mapsto n^j (r_m^\ell)^{\theta n}, \theta \ne 0,$ extends to a well-defined surjective map from
the set of equivalence classes of type
(\ref{eq:form1}) and (\ref{eq:form2}) to the set $U_Q \cup W_R$, proving that $|U_Q \cup W_R|\le M$.
Hence the C-finite degree of $f(n)$ is at most $M$ in all cases, and
the remaining assertions follow
immmediately.
$\Box$

It is possible to give another bound on the C-finite degree that is sometimes sharper than the one in Corollary
\ref{gendegbound}.

\begin{corollary}\label{gendegboundx}
Suppose that $F_1, F_2,\dots, F_k$ are as in Theorem \ref{th:mpl}
where each $F_i$ is C-finite of degree $D(i)$, with $d(i)$
distinct roots.
Then the C-finite degree of $f(n)$ is bounded by
\begin{equation}\label{gendegbound1}
\left(\prod_{\ell\in Q} d(\ell)\right)  \left( 1+ \sum_{\ell\in Q}
\left(\frac{D(\ell)}{d(\ell)}-1\right) \right)+
\left(\prod_{\ell\in R} d(\ell) \right) \left( 1+ \sum_{\ell\in R}
\left(\frac{D(\ell)}{d(\ell)}-1\right) \right),
\end{equation}
where $Q=\{\ell\;|\;a_\ell+b_\ell\not=0\}$ and
$R=\{\ell\;|\;a_\ell\not=0\}$. This last expression is in turn
bounded by
\begin{equation}\label{gendegbound2}
\prod_{\ell\in Q} D(\ell) + \prod_{\ell\in R} D(\ell).
\end{equation}
In all of these expressions, the empty product is taken to be equal to
$1$.
\end{corollary}
\noindent
{\bf Proof:}
Again, we first consider the case where no
product of the form $\Theta = \prod_{\ell=1}^k (r_{m_\ell}^{(\ell)})^{b_\ell}$ equals $1$.
To establish the bound in (\ref{gendegbound1}), note that
\begin{eqnarray}\label{ineq32}
|U_Q| \leq
\sum_{(m_\ell)_{\ell\in Q}} \left(1+ \sum_{\ell\in Q} (e_{m_\ell}^{(\ell)}-1)\right) &=&
\left(\prod_{\ell\in Q} d(\ell) \right)+ \left(\sum_{\ell\in Q}\bigg( D(\ell)\prod_{j\ne \ell,j\in Q} d(j)-
\prod_{\ell\in Q} d(\ell)\bigg)\right)\\
&=&\left(\prod_{\ell\in Q} d(\ell)\right)  \left( 1+\sum_{\ell\in Q} \bigg(\frac{D(\ell)}{d(\ell)}-1\bigg) \right)\nonumber
\end{eqnarray}
where the first summation is over all sequences $(m_\ell)_{\ell \in Q}$, satisfying
$1 \le m_\ell \le d(\ell)$ for $\ell\in Q$. The second term in (\ref{gendegbound1})
is handled similarly, and the result follows.
To obtain the bound in (\ref{gendegbound2}), it suffices to show that
$$
\left(\prod_{\ell\in Q} d(\ell)\right)  \left( 1+ \sum_{\ell\in Q}
\left(\frac{D(\ell)}{d(\ell)}-1\right) \right) \le
\prod_{\ell\in Q} D(\ell)
$$
To see this,
divide both sides by $\prod_{\ell\in Q} d(\ell)$, and write $D(\ell)/d(\ell) = 1+t_\ell$, where
$t_\ell \ge 0$. It remains to show that
$
1+ \sum_{\ell\in Q} t_\ell \le \prod_{\ell \in Q}(1+t_\ell)
$
for nonnegative $t_\ell$, which is obvious.
This completes the
proof of Corollary \ref{gendegbound} in the first case.

Next suppose that there are products of the form
$\Theta = \prod_{\ell=1}^k (r_{m_\ell}^{(\ell)})^{b_\ell}$ equal to $1$.
Let $W_R$ be defined as in the proof of Corollary \ref{gendegbound}.
We will show that $|U_Q\cup W_R|$ is bounded by (\ref{gendegbound1}).
Arguing as in the first case, we obtain that $|U_Q|$ is less than or equal
to the first summand in (\ref{gendegbound1}). Before computing the $W_R$
contribution, it is convenient to define
$$
W_R^0 \;=\; \{ n^h \prod_{\ell\in R} (r_{m_\ell}^{(\ell)})^{a_\ell n} \;|\; H_0 < h \le H\}
$$
where $H_0 = 1 + \sum_{\ell\in R} (e_{m_\ell}^{(\ell)} - 1)$. We claim that
$$
|U_Q \cap W_R| \;\ge\; |W_R^0|.
$$
To see this, suppose that $n^h \prod_{\ell\in R}
(r_{m_\ell}^{(\ell)})^{a_\ell n}\in W_R^0$. Since $h>H_0$, there
must exist indices $m^*_\ell$ with $\ell\not\in R$ such that
$\prod_{\ell\in R} (r_{m_\ell}^{(\ell)})^{b_\ell}
\prod_{\ell\not\in R} (r_{m^*_\ell}^{(\ell)})^{b_\ell}= 1$. Assume that these have
been chosen so that the right hand side of (\ref{h0def}) is maximized.
Then
$$
n^{h-H_0} \prod_{\ell\in R} (r_{m_\ell}^{(\ell)})^{a_\ell n}
=
n^{h-H_0} \prod_{\ell\in R} (r_{m_\ell}^{(\ell)})^{(a_\ell+b_\ell) n}
\prod_{\ell\not\in R} (r_{m^*_\ell}^{(\ell)})^{b_\ell n}
$$
is an element of $U_Q\cap W_R$, and it is clear that this map is injective. It follows that
\begin{eqnarray*}
|U_Q \cup W_R| &=& |U_Q| + |W_R| - |U_Q \cap W_R| \\
&\le&|U_Q| + |W_R| - |W_R^0| \\
&=&|U_Q| + |W_R-W_R^0|.
\end{eqnarray*}
Finally, we have
$$
|W_R-W_R^0| \le \sum_{(m_\ell)_{\ell\in R}} (1+\sum_{\ell\in R} (e_{m_\ell}^{(\ell)}-1)),
$$
which is less than or equal to the second summand in (\ref{gendegbound1}), by the argument presented
in the first case.  The bound (\ref{gendegbound2}) follows as before, and the proof is complete.
$\Box$

We remark that the second case above can also be derived from the first case by a continuity argument.

\section{A striking property of indefinite summation}

When $a_1=a_2=\cdots=a_k=0$ in (\ref{eq:gensum}) we are doing indefinite summation, i.e.,
the problem is equivalent to finding a function $S(n)$ (called an \textit{indefinite sum function})
such that
\begin{equation}\label{indef}
S(n)-S(n-1) \;=\; F_1(b_1n +c_1) \cdots F_k(b_kn +c_k)
\end{equation}
In this section we will show that, if no products of the form
$\prod_{i=1}^k(r_{m_i}^{(i)})^{b_i}$
are equal to $1$,
there exists an indefinite sum function $S(n)$ expressible (formally) as a linear
combination of monomials in the $F_i$ with coefficients that are independent of the initial conditions
satisfied by the various $F_i$.

Theorem \ref{inits} proved this result (which we will call the
{\it independence property}) for power sums $\sum_j F(j)^p$,
assuming that the roots $r_i$ are distinct, and it is easy to
extend that proof to the case of summands of the form $F_1(b_1n
+c_1) \cdots F_k(b_kn +c_k)$, with different $F_i$,
as long as the roots are distinct. Among other things, this
section extends Theorem \ref{inits}
to the multiple root case, where the proof
turns out to be considerably more difficult.

Similar questions have been considered in
\cite{ru1}, where the author describes a method (different from ours) that
finds indefinite sum functions
in some but not all cases\footnote{The authors thank the referee for bringing this reference to their
attention, and also for several other helpful suggestions and remarks.},
but does not delineate
these cases with a theorem.
In fact, the method developed in
\cite{ru1} assumes the independence property;
the author asserts (in somewhat vague terms, and without proof) that it holds if no products
$\prod_i (r_{m_i}^{(i)})^{b_i}$ are equal to $1$.

Theorem \ref{bigcor} (below) establishes this result, and in full generality.
It dispenses with the assumption of distinct roots, and makes appropriate modifications
in the case where products of the form $\prod_i (r_{m_i}^{(i)})^{b_i}$ are equal to $1$.
The first step in the proof
is to identify a collection of functions analogous to the terms
$\lambda_mr_m^n, 1\leq m\leq d$ that appear in the proof of Theorem \ref{inits}, which
can be expressed as linear combinations of the $F(n+i)$ with coefficients that do not depend
on the initial conditions.

\begin{lemma}\label{indepphi}
Suppose that
$$
F(n) = \sum_{m=1}^{d} \sum_{h=0}^{e_{m}-1}\,\lambda_{m,h}\,n^h \,(r_m)^n,
$$
where the $r_m$ are distinct and nonzero, and
$\lambda_{m,e_m-1}\not=0$ for all $m$. For $m=1,\dots,d$ and $h=0,\dots, e_m-1$, define
$$
\Phi_h^{m}(n) \;=\;  \sum_{i=h}^{e_m-1}
{i\choose h} \lambda_{m,i} n^{i-h} (r_{m})^n
$$
Let $D=\sum_m e_m$.
Then each of the functions $\Phi^m_h(n)$ may be expressed as a linear combination
of the functions $F(n), F(n+1), \dots, F(n+D-1)$, with coefficients that do not depend on the
values of $\lambda_{m,h}$, i.e.,
those coefficients do not depend on the initial conditions of $F$.
\end{lemma}
\noindent{\bf Proof: } Let ${\bf F}$ be the column vector of
length $D$ whose $i$th element is equal to $F(n+i),
0\leq i \leq D-1$. Let ${\bf R}$ and ${\bf \Phi}$ be column vectors of length
$D$, each indexed by $(m,h), 1\leq m \leq d, 0\leq h \leq e_m-1$, where the
$(m,h)$th element of ${\bf R}$ is equal to $n^h \,(r_m)^n$ and the
$(m,h)$th element of ${\bf \Phi}$ is equal to $\Phi_{e_m-1-h}^m(n)$. Then we have
$$
{\bf F} \;=\; {\bf M}\,{\bf R}
$$
where ${\bf M}$ is a matrix with rows indexed by $0,\dots,D-1$ and columns
indexed by $(m,h)$, whose entry in row $t$ and column $(m,h)$ is equal to
$$
(r_m)^t \sum_{i=h}^{e_m-1} {i \choose h} \lambda_{m,i} t^{i-h}
$$
We have the factorization
$$
{\bf M} \;=\; {\bf M_0}\,{\bf \Lambda}
$$
where ${\bf M_0}$ is a matrix indexed as in ${\bf M}$, whose entry in row $t$ and
column $(m,h)$ is
$$
(r_m)^t \, t^{e_m-1-h},
$$
and ${\bf \Lambda}$ is a matrix with both rows and columns indexed by $(m,h)$, whose
entry is row $(m,h)$ and column $(m',h')$ is equal to zero if $m\not=m'$ or $h'>h$,
and is otherwise equal to
$$
\lambda_{m,e_m-1-(h-h')}{e_m-1-(h-h')\choose h'}.
$$
One can also check that
$$
{\bf \Phi} \;=\; {\bf \Lambda} \, {\bf R}.
$$
We note that ${\bf M_0}$ does not depend on the values of
$\lambda_{m,h}$. Furthermore, after factoring out various powers of the roots $r_m$, it is
column-equivalent to the transpose of the matrix whose determinant was computed in
Theorem \ref{nicematrix}. Hence,
since we are assuming the $r_m$ to be distinct and nonzero, ${\bf M_0}$ is
nonsingular. Thus we can write
$$
{\bf M_0}^{-1}\,{\bf F} \;=\; {\bf \Lambda}\,{\bf R}\;=\; {\bf \Phi},
$$
and the lemma follows.
$\Box$.

As an illustration of the previous lemma and its proof, consider the function
$$
F(n) = a r_1^n + b n r_1^n + c n^2 r_1^n + d r_2^n + e n r_2^n
$$
Then
\begin{eqnarray*}
\Phi^{1}_0(n) &=& a \, r_1^n + b n\, r_1^n  + cn^2\, r_1^n\\
\Phi^{1}_1(n) &=&b \, r_1^n + 2c n\, r_1^n \\
\Phi^{1}_2(n) &=& c\, r_1^n\\
\Phi^{2}_0(n) &=& d\, r_2^n  + e n \,r_2^n \\
\Phi^{2}_1(n) &=& e\,r_2^n.
\end{eqnarray*}
We have
$$
{\bf M\,R} \;=\;
\left(
\begin{array}{ccccc}
a & b & c & d & e\\
(a+b+c)r_1 & (b+2c)r_1 & cr_1 & (d+e)r_2 & er_2 \\
(a+2b+4c)r_1^2 & (b+4c)r_1^2 & cr_1^2 & (d+2e)r_2^2 & e r_2^2\\
(a+3b+9c)r_1^3 & (b+6c)r_1^3 & c r_1^3& (d+3e)r_2^3 & e r_2^3\\
(a+4b+16c)r_1^4 & (b+8c)r_1^4 & c r_1^4& (d+4e) r_2^4& er_2^4
\end{array}
\right)
\left(
\begin{array}{c}
r_1^n\\
n r_1^n\\
n^2 r_1^n\\
r_2^n\\
n r_2^n
\end{array}
\right)
\;=\;
\left(
\begin{array}{c}
F(n)\\
F(n+1)\\
F(n+2)\\
F(n+3)\\
F(n+4)
\end{array}
\right)
$$
and
$$
{\bf M}\;=\;
{\bf M_0}\,{\bf \Lambda} \;=\;
\left(
\begin{array}{ccccc}
0 & 0 & 1 & 0 & 1\\
 r_1& r_1 & r_1 & r_2 & r_2\\
4r_1^2 & 2r_1^2 & r_1^2 & 2r_2^2 & r_2^2\\
9r_1^3 & 3r_1^3 & r_1^3& 3r_2^3 & r_2^3\\
16 r_1^4& 4r_1^4 & r_1^4 & 4r_2^4 & r_2^4
\end{array}
\right)
\left(
\begin{array}{ccccc}
c & 0 & 0 & 0 & 0\\
b & 2c & 0 & 0 & 0\\
a & b & c & 0 & 0\\
0 & 0 & 0 & e & 0\\
0 & 0 & 0 & d & e
\end{array}
\right)
$$
Also
$$
{\bf \Phi}\;=\;
\left(
\begin{array}{c}
\Phi^1_2(n)\\
\Phi^1_1(n)\\
\Phi^1_0(n)\\
\Phi^2_1(n)\\
\Phi^1_0(n)
\end{array}
\right)\;=\;
\left(
\begin{array}{ccccc}
c & 0 & 0 & 0 & 0\\
b & 2c & 0 & 0 & 0\\
a & b & c & 0 & 0\\
0 & 0 & 0 & e & 0\\
0 & 0 & 0 & d & e
\end{array}
\right)
\left(
\begin{array}{c}
r_1^n\\
n r_1^n\\
n^2 r_1^n\\
r_2^n\\
n r_2^n
\end{array}
\right)
\;=\; {\bf \Lambda}\,{\bf R}
$$
Thus the relation ${\bf \Phi}={\bf M_0}^{-1} {\bf F}$ expresses each $\Phi_h^m(n)$ as a
linear combination of $F(n),\dots,F(n+4)$, with coefficients that do not depend on $a,b,c,d$,
and $e$.

To prove our main result, we also need the following Lemma.
\begin{lemma}\label{sumformula}
If $p\geq 0$, define $S_p(n,x) = \sum_{j=0}^{n-1} j^p\,x^j$. If $x\not=1$, then
\begin{eqnarray*}
S_p(n,x)&=& B_p(x) - x^n \big( \sum_{k=0}^p {p \choose k} n^k B_{p-k}(x) \big)\\
&=& B_p(x) - x^n \big( (n+Q)^p \bigg|_{Q^k\to B_k(x)} \big)
% &\eqdef& H_p(0,x)- H_p(n,x)
\end{eqnarray*}
where $B_k(x) = A_k(x) / (1-x)^{k+1}$, and $A_k(x)$ denotes the
Eulerian polynomial of degree $k$. Here, the notation $Q^k\to
B_k(x)$ means ``replace $Q^k$ by $B_k(x)$ throughout".
If $x=1$, then $S_p(n,x)$ is a (well-known) polynomial
of degree $p+1$.
\end{lemma}

\noindent{\bf Proof: } The identity
$$
\sum_{j=0}^\infty j^p x^j = \frac{A_p(x)}{(1-x)^{p+1}}
$$
is classical (see e.g. \cite{co}), and immediately
implies the relation $xB_p'(x) = B_{p+1}(x)$, as a formal power series identity.
The lemma now follows easily by induction. $\Box$

The next Lemma provides an explicit form for our indefinite sums,
expressing the result in terms of the functions $\Phi_h^m(n)$ defined
in Lemma \ref{indepphi}.
\begin{lemma}\label{bigformula}
Suppose that for each $i=1,\dots,k$,
$$
F_i(n) = \sum_{m=1}^{d(i)} \sum_{h=0}^{e^{(i)}_{m}-1}\,\lambda^{(i)}_{m,h}\,n^h \,(r_m^{(i)})^n
$$
with $\lambda^{(i)}_{m,e^{(i)}_m-1}\not=0$ for all $i$ and $m$. For $1 \le i \le k$, let
$\delta_{i,m}= e^{(i)}_m-1$. Then
\begin{equation}\label{eqsum}
\sum_{j=0}^{n-1} F_1(j)F_2(j)\cdots F_k(j) \;=\; \Psi(n) + S(0) - S(n)
\end{equation}
where $\Psi(n)$ is a polynomial, and
\begin{equation}\label{eqk}
S(n) \;=\; \sum_{r_{m_1}^{(1)}\cdots r_{m_k}^{(k)}\ne 1} \sum_{s=0}^{\delta_{1,m_1}+\cdots+\delta_{k,m_k}}
B_s(r_{m_1}^{(1)}\dots r_{m_k}^{(k)}) \left\{ \sum_{t_1+\cdots+t_k=s}
\Phi_{t_1}^{1,m_1}(n)\cdots \Phi_{t_k}^{k,m_k}(n) \right\}.
\end{equation}
If no product of the form $r_{m_1}^{(1)}\cdots r_{m_k}^{(k)}$
equals $1$, then $\Psi(n)\equiv 0$; otherwise $\Psi(n)$ has degree at most
$1 + \sum_i\Delta(i)$, where
$\Delta(i)= \max_{1\le j\le d(i)}(e_j^{(i)}-1)$, as defined in Theorem \ref{th:mpl}. In formula
(\ref{eqk}),
$B_s(x) = A_s(x)/(1-x)^{s+1}$, as defined
in Lemma \ref{sumformula},  and
for $\ell=1,\dots,k$,
$$
\Phi_t^{\ell,m_\ell}(n) \;=\;  \sum_{i=t}^{\delta_{\ell,m_\ell}} {i\choose t}
\lambda^{(\ell)}_{m,i} n^{i-t} (r^{(\ell)}_{m_\ell})^n,
$$
as defined in Lemma \ref{indepphi}.
\end{lemma}
Note that, in expression (\ref{eqk}), the coefficients $B_s(r_{m_1}^{(1)}\dots r_{m_k}^{(k)})$ depend
only on the roots $r^{(i)}_m$, and hence by Lemma \ref{indepphi}, $S(n)$ is
independent of the initial conditions of the $F_i$.

\noindent{\bf Proof: }
To simplify notation, we will first consider the
case where the $F_i$ each have a single root, i.e.
\begin{eqnarray*}
F_i(n) &=& \big( \lambda_0^{(i)} + \lambda_1^{(i)} n  + \cdots + \lambda_{\delta_i}^{(i)}
n^{\delta_i}\big)(r^{(i)})^n
\end{eqnarray*}
where $r^{(i)}$ is a root of multiplicity $e^{(i)}=\delta_i+1$. We have
\begin{eqnarray*}
\sum_{j=0}^{n-1} F_1(j) \cdots F_k(j) &=&
\sum_{j=0}^{n-1} \,\big( \sum_{i=0}^{\delta_1} \lambda_i^{(1)} j^i (r^{(1)})^j \big)
\cdots
\big( \sum_{i=0}^{\delta_k} \lambda_i^{(k)} j^i (r^{(k)})^j \big)\\
&=&\sum_{j=0}^{n-1} \;\sum_{s=0}^{\delta_1+\cdots+\delta_k} \sum_{i_1+\cdots +i_k=s}
\lambda_{i_1}^{(1)}\cdots \lambda_{i_k}^{(k)}j^s (r^{(1)}\cdots r^{(k)})^j\\
&=&\sum_{s=0}^{\delta_1+\cdots+\delta_k} \sum_{i_1+\cdots +i_k=s}
\lambda_{i_1}^{(1)}\cdots \lambda_{i_k}^{(k)}\,\sum_{j=0}^{n-1}\, j^s (r^{(1)}\cdots r^{(k)})^j.
\end{eqnarray*}
Now for any integer $p\geq 0$, define
\begin{equation}\label{hdef}
H_p(n,x) \;=\; x^n \sum_{k=0}^p {p \choose k} n^k B_{p-k}(x),
\end{equation}
where $B_k(x)$ is as defined in Lemma \ref{sumformula}. Thus, by that Lemma, we have
$$
\sum_{j=0}^{n-1} j^p x^j \;=\; H_p(0,x) - H_p(n,x)= x^n (n+Q)^p \bigg|_{Q^k\to B_k(x)},
$$
for any $x\ne 1$. Continuing with the above calculation, we obtain
$$
\sum_{j=0}^{n-1} F_1(j)\cdots F_k(j) \;=\;
\sum_{s=0}^{\delta_1+\cdots+\delta_k} \sum_{i_1+\cdots +i_k=s}
\lambda_{i_1}^{(1)}\cdots \lambda_{i_k}^{(k)}\,\left(H_s(0, r^{(1)}\cdots r^{(k)}) -H_s(n, r^{(1)}\cdots r^{(k)})\right)
$$
provided that $r^{(1)}\cdots r^{(k)}\ne 1$; otherwise the sum is equal to $\Psi(n)$, a polynomial of degree at most equal
to $1 + \sum_i \delta_i = 1 + \sum_i( e^{i}-1)$.
In the former case, write
$$
S(n) \;=\;\sum_{s=0}^{\delta_1+\cdots+\delta_k} \sum_{i_1+\cdots +i_k=s}
\lambda_{i_1}^{(1)}\cdots \lambda_{i_k}^{(k)}\,H_s(n, r^{(1)}\cdots r^{(k)})
$$
so that
\begin{equation}\label{sumcases}
\sum_{j=0}^{n-1} F_1(j)\cdots F_k(j) \;=\;
\begin{cases} S(0)-S(n) & \text{if } r^{(1)}\cdots r^{(k)}\ne 1, \text{ and}\\
\Psi(n) & \text{otherwise}.
\end{cases}
\end{equation}
Using the alternate form in (\ref{hdef}), we can write
\begin{eqnarray*}
S(n) &=&\sum_{s=0}^{\delta_1+\cdots+\delta_k} \sum_{i_1+\cdots +i_k=s}
\lambda_{i_1}^{(1)}\cdots \lambda_{i_k}^{(k)}  (r^{(1)}\cdots r^{(k)})^n(n+Q)^s  \;\bigg|_{Q^s \to B_s(r^{(1)}\cdots r^{(k)})}\\
&=& \bigg( \sum_{i=0}^{\delta_1} \lambda_{i}^{(1)} (n+Q)^i\bigg) (r^{(1)})^n \cdots
\bigg( \sum_{i=0}^{\delta_k} \lambda_{i}^{(k)} (n+Q)^i\bigg) (r^{(k)})^n \;
\;\bigg|_{Q^s \to B_s(r^{(1)}\cdots r^{(k)})}\\
&=&\bigg( \sum_{i=0}^{\delta_1} \lambda_{i}^{(1)} \sum_{u=0}^i {i \choose u} n^{i-u} Q^{u}\bigg) (r^{(1)})^n\cdots
\bigg( \sum_{i=0}^{\delta_k} \lambda_{i}^{(k)}  \sum_{u=0}^i {i \choose u} n^{i-u} Q^{u}   \bigg) (r^{(k)})^n \;
\;\bigg|_{Q^s \to B_s(r^{(1)}\cdots r^{(k)})}\\
&=&\bigg( \sum_{u=0}^{\delta_1}  \sum_{i=u}^{\delta_1} {i \choose u}\lambda_{i}^{(1)} n^{i-u}(r^{(1)})^n Q^{u}\bigg) \cdots
\bigg( \sum_{u=0}^{\delta_k}   \sum_{i=u}^{\delta_k} {i \choose u}\lambda_{i}^{(k)} n^{i-u} (r^{(k)})^n Q^{u}   \bigg)  \;
\;\bigg|_{Q^s \to B_s(r^{(1)}\cdots r^{(k)})}\\
&=&
\bigg( \sum_{u=0}^{\delta_1} \Phi^1_u(n) Q^u \bigg) \cdots \bigg(\sum_{v=0}^{\delta_k} \Phi^k_v(n) Q^v\bigg)
 \;\;\bigg|_{Q^s \to B_s(r^{(1)}\cdots r^{(k)})}
\end{eqnarray*}
Combining this last expression with (\ref{sumcases}) yields the statement of Lemma \ref{bigformula} in the case where each
recurrence has a single root.
For the general case (i.e., when there are
multiple roots for each $F_i$), one can collect terms in the expansion of $F_1(j)F_2(j)\cdots F_k(j)$
corresponding to each choice $r_{m_1},r_{m_2},\dots,r_{m_k}$ of a sequence of roots from each $F_i$, and
an expression for $K(n)$ of the form (\ref{eqsum}) results, with $S(n)$ as in
(\ref{eqk}). This completes the proof of Lemma  \ref{bigformula}.
$\Box$.

Combining Lemmas \ref{indepphi} and \ref{bigformula}, we obtain the following Theorem,
which is the main result of this section.
\begin{theorem}\label{bigcor}
If $F_1(n), F_2(n),\dots,F_k(n)$ satisfy the hypotheses of Lemma \ref{bigformula}, then the sum
$$
\sum_{j=0}^{n-1} F_1(b_1 j+ c_1)F_2(b_2 j+ c_2)\cdots F_k(b_k j + c_k)
$$
may be expressed as
$\Psi(n) -S(n)$, where $\Psi(n)$ is a polynomial and $S(n)$ is a linear combination
of monomials of the form
$$F_1(b_1n+i_1)F_2(b_2n+i_2)\cdots F(b_k n+i_k) \qquad (0\leq i_\nu \leq D(\nu)-1;1\le \nu \le k)$$ with
$D(\nu) = \sum_m e^{(\nu)}_m$, such that the coefficients of that linear combination
are independent of the initial conditions of the $F_i$.
If no product $(r_{m_1}^{(1)})^{b_1}\cdots (r_{m_k}^{(k)})^{b_k}$ of the associated roots equals $1$,
then $\Psi(n)$ is a constant; otherwise
it has degree at most $1 + \sum_i \Delta(i)$.
In general, $\Psi(n)$ will depend on the
initial conditions of the $F_i$.
\end{theorem}

\section{More examples}
We will give some examples illustrating the results in the previous two sections.
\subsection{A mixed convolution}
Let $F(n)$ denote the $n$th Fibonacci number, and let $G(n)$ be defined
by the subword-avoiding recurrence (\ref{eq:subwrecur}) with $A=3$, in
other words $G(0)=1, G(1)=3, G(2)= 9$, and $G(n)= 3G(n-1)-G(n-3)$ for
$n>2$. Then we have the following identity:
$$\sum_{j=0}^n j\, F(j) G(n-j)\;=\;
18 G(n+1) - (9G(n)+ 5 G(n+2) + 3 F(n) +n F(n) + n F(n+1))
$$
The target monomials in this case are
$$
F(n),\, n F(n),\, F(n+1),\, n F(n+1), \,G(n), \,G(n+1), \,G(n+2)
$$
and the (unique) solution is obtained by solving a system of $7$ equations
in $7$ unknowns. Here we are applying Theorem \ref{th:mpl}
in the case where no product of the form $\prod r_m^b$ equals $1$, with $F_1(n)=n,F_2(n)=F(n)$, and
$F_3(n)=G(n)$. The sum is C-finite of degree $7$, and this
degree achieves the bounds given in Corollaries
\ref{gendegbound} and \ref{gendegboundx}.

\subsection{The independence property}
We will give two examples of indefinite summations illustrating the phenomena described in
Theorems \ref{inits} and \ref{bigcor}.
Consider the sum
$ \sum_{j=0}^{n-1} F(j)^3
$
where $F(n)$ satisfies the Fibonacci recurrence with initial conditions
$F(0)=p, F(1)=1$. The target monomials are
$$
1,\,F(n)^3,\, F(n)^2 F(n+1),\,F(n) F(n+1)^2,\,F(n+1)^3
$$
and we obtain the identity
$$
\sum_{j=0}^{n-1} F(j)^3\;=\;
\frac{1}{2}(1-3p +3p^3) - \frac{3}{2}F(n)^3 + \frac{3}{2} F(n)F(n+1)^2 - \frac{1}{2} F(n+1)^3,
$$
which is of the form $S(0)-S(n)$ in the notation of Theorem \ref{bigcor}. Since the associated
roots are distinct and no product equals $1$, this case is covered by Theorem \ref{inits}.

Next consider the sum $\sum_{j=0}^{n-1} F(j)^4$, for which we now have products of roots equal to $1$. The
target monomials are
$$
1,\,n\,F(n)^4,\, F(n)^3 F(n+1),\,F(n)^2 F(n+1)^2,\,F(n)F(n+1)^3, F(n+1)^4
$$
and we obtain the general solution
\begin{eqnarray*}
&\frac{1}{25}\big(A_0 \,+\, A_1 n \,+\, 52F(n)^3 F(n+1)-22F(n)^2 F(n+1)^2-36F(n)F(n+1)^3+19 F(n+1)^4\big)&\\
&\qquad + K ((-1+p+p^2)^2 - (F(n+1)^2 - F(n)^2 - F(n)F(n+1))^2&
\end{eqnarray*}
where
$$
A_0= -19 + 36p + 22 p^2 -52 p^3, \quad A_1 = 6(-1+p+p^2)^2,
$$
and $K$ is an arbitrary constant.
If $p^2 + p -1 \ne 0$, the polynomial term $\Psi(n)$ in Theorem \ref{bigcor} will
always have degree $1$.

\subsection{Partial summation of series}
Consider the sum
$$ \sum_{j=0}^{n-1} F(j) x^j
$$
where $F(n)$ is the $n$th Fibonacci number and $x$ is an
indeterminate. The summand is a product
of two C-finite sequences, one of degree two and the other of degree one.
Following Theorem \ref{th:mpl}, we construct a list of target
of target monomials $1, F(n)x^n$, and $F(n+1) x^n$, and from these we
obtain the identity
$$
\sum_{j=0}^{n-1}F(j)x^j=\frac{x}{1-x-x^2}-x^n \bigg(\frac{1-x}{1-x-x^2}F(n)+
\frac{x}{1-x-x^2}F(n+1)\bigg).
$$
This identity quantifies the remainder term in the Fibonacci generating function (an
equivalent result appears as problem 1.2.8.21 in
\cite{dek}).
Our approach can be easily extended; for example, using
$1, F(n)^2 x^n, F(n)F(n+1) x^n$, and $F(n+1)^2 x^n$ as target monomials and solving
four equations in four unknowns, we
obtain the partial summation formula
\begin{equation}\label{eq:f2gf}
\sum_{j=0}^{n-1} F(j)^2 x^j = \frac{x(1-x)}{1-2x-2x^2+x^3} -
x^n \,R_n(x)
\end{equation}
where
\begin{equation}\label{eq:f2gf2}
R_n(x)= \frac{(1-2x-x^2)F(n)^2 + 2x^2F(n)F(n+1) + x(1-x)F(n+1)^2}
{1-2x-2x^2+x^3}
\end{equation}
The first term in (\ref{eq:f2gf}) is the full generating function for squares
of Fibonacci numbers. A formula for the full generating function
for all powers $p$ appears in \cite{r} (see also \cite{dek}, problem 1.2.8.30).

We note that, to obtain
(\ref{eq:f2gf}) and (\ref{eq:f2gf2}) by this method, it was only necessary
to know the first four values of the sum, and also that
$F$ satisfies some $2$-term recurrence with constant coefficients.

\subsection{The degree of the polynomial multiplier}

In Theorem \ref{th:mpl} we gave a set of monomials in terms of which the sum can be expressed,
in the general case of repeated roots. In those monomials a polynomial factor $\psi_{i_1\dots i_k}$
appears, and the degree of that polynomial was found to be at most $1+\sum \{\Delta(i):a_i=0\}$. We remark
here that the $C$-finite function $F(n)=n^p$, for positive integer $p$, shows that this upper bound
can be achieved. For here we have
\[k=1, d(1)=1,\,r_1^{(1)}=1,\,e_1^{(1)}=p+1,\,\Delta(1)=p,\,D(1)=p+1,\,(a_1,b_1,c_1)=(0,1,0).\]
The monomials in the list (\ref{eq:form1}) are all of degree $p$, and those in the list
(\ref{eq:form2}) are of degree equal to the degree of $\psi_1(n)$.
The maximum allowable degree of the latter is $1+\Delta(1)=p+1$, and
in this case $\psi$ is of degree $p+1$ since the sum obviously is.

\subsection{An example from the theory of partitions}
Our algorithm for summation of $C$-finite sequences can sometimes
have by-products that are more interesting than the particular
problem being solved. A small example of this is given here.
Suppose $p_5(n)$ is the number of partitions of $n$ into $\le 5$
parts. Then $\sum_{n\ge 0}p_5(n)z^n=\left((1-z)(1-z^2)(1-z^3)(1-z^4)(1-z^5)\right)^{-1}$,
so $p_5$ is $C$-finite of degree $\le 15$.
We asked our Mathematica program \cite{gwmma} to find $f(n)=\sum_{0\le j\le n-1}p_5(j)$.
In addition to giving the answer, a number of the arbitrary constants that
are used to form linear combinations with target monomials were left unassigned,
and since the coefficient of every such unassigned constant must of course vanish,
one has found an identity. On this occasion we chose one symmetrical looking such identity from the output, namely
\begin{eqnarray*}
&& p_5(n-10 ) + 4\,p_5(n-9 ) + 9\,p_5(n-8) + 15\,p_5(n-7) +
  20\,p_5(n-6) + 22\,p_5(n-5) + 20\,p_5(n-4)\\
&&\qquad\qquad   + 15\,p_5(n-3) + 9\,p_5(n-2) + 4\,p_5(n-1) + p_5(n)={n+4\choose 4}.
\end{eqnarray*}
The coefficients on the left side are recognized as the numbers of permutations of 5 letters
that have $k$ inversions. From this one might suspect that we have generally,
\begin{equation}\label{eq:ident}
\sum_jb(k,j)p_k(n-j)={n+k-1\choose n},
\end{equation}
where $b(k,j)$ is the number of $k$-permutations that have exactly $j$ inversions and $p_k(m)$
is the number of partitions of $m$ into parts $\le k$. A proof of this identity by generating
functions is quite trivial. Here is a bijective proof, that is, a bijection between pairs
consisting of a permutation of $k$ letters with $j$ inversions and a partition of $n-j$
into $\le k$ parts, on the one hand, and one of the ${n+k-1\choose n}$ compositions
of $n$ into $k$ nonnegative parts, on the other. Take such a composition $X$ of $n$
into $k$ parts. Perform a ``modified bubble sort,"
whereby whenever one sees an adjacent pair $x\,y$ with $x<y$, it is replaced by $(y-1)x$.
Keep doing this until there are no adjacent pairs $x<y$, i.e., until a partition $\lambda$ (perhaps followed by 0's)
is obtained.  Call the resulting permutation of positions $\sigma$.
Then $\lambda$ and $\sigma$ are uniquely determined by $X$, and the correspondence
is bijective. $\Box$

We do not claim novelty for this result or its proof, but offer
it only as an example of the usefulness that our algorithms can have in the discovery process.


\begin{thebibliography}{aaa}
\bibitem{co} L. Comtet, {\it Advanced Combinatorics}, D. Reidel, Dordrecht, 1974.
\bibitem{fh} R. P. Flowe, G. A. Harris, A note on generalized Vandermonde
determinants, SIAM J. Matrix Anal. Appl. {\bf 14 4} (1993), 1146-1151.
\bibitem{jordan} Charles Jordan, {\it Calculus of Finite Differences},
Chelsea, New York, 1950.
\bibitem{gwmma} Curtis Greene, Herbert S. Wilf, \texttt{CFSum.nb}, (Mathematica notebook),
\texttt{<http://www.haverford.edu/math/cgreene/cfsum.nb>},
\texttt{<http://www.math.upenn.edu/}$\sim$\texttt{wilf/website/cfsum.nb>}.
\bibitem{dek}
Donald E. Knuth, {\it The Art of Computer Programming}, Addison-Wesley, Reading,
MA, 1969, Vol. 1, p. 84 (exercises 1.2.8.21 and 1.2.8.30), and p. 491, p. 492 (solutions).
\bibitem{k} Christian Krattenthaler, Advanced Determinant Calculus,
S\'eminaire Lotharingien Combin. \textbf{42} (``The Andrews Festschrift'')
(1999), Article B42q, 67 pp.
\bibitem{pwz} Marko Petkov\v sek, Herbert S. Wilf, and Doron
Zeilberger,
$A=B$, A K Peters Ltd., Wellesley, MA, 1996.
\bibitem{r}J. Riordan, Generating functions for powers of Fibonacci numbers,
Duke. Math. J. 29 (1962) 5-12.
\bibitem{ru1} David L. Russell, Sums of products of terms from linear recurrence
sequences, Discrete Math \textbf{28} (1979), 65-79.
\bibitem{z1} Doron Zeilberger, A holonomic systems approach to special
functions identities, J. Comput. Appl. Math. \textbf{32} (1990), no. 3,
321--368.
\end{thebibliography}
\end{document}